\documentclass[aap]{imsart}

\RequirePackage{amsthm,amsmath,amsfonts,amssymb,bbm,color,tikz}
\RequirePackage[numbers]{natbib}
\RequirePackage[colorlinks,citecolor=blue,urlcolor=blue]{hyperref}

\startlocaldefs

\allowdisplaybreaks

\newtheorem{theorem}{Theorem}[section]

\newtheorem{lemma}[theorem]{Lemma}
\newtheorem{proposition}[theorem]{Proposition}

\theoremstyle{definition}

\newtheorem{definition}[theorem]{Definition}

\newtheorem{remark}[theorem]{Remark}

\let\plainqed\qedsymbol
\newcommand{\claimqed}{$\lrcorner$}

\newcommand{\Expect}[1]{\operatorname{\mathbb{E}}\left[#1\right]}
\newcommand{\Identity}{\operatorname{\mathbb{I}}}

\newcommand{{\LPC}}{\textbf{LPC}}
\newcommand{\pr}[1]{\operatorname{\mathbb{P}}\left[#1\right]}

\newcommand{\LL}{\mathcal{L}}

\newcommand{\G}{\mathcal{G}}

\newcommand{\DPA}{\operatorname{DPA}}

\newcommand{\PP}{\mathcal{P}}
\newcommand{\II}{\mathcal{I}}
\newcommand{\B}{\mathcal{B}}
\newcommand{\TT}{\mathbb{T}}
\newcommand{\oa}{\mathbf{i}}
\newcommand{\equald}{\stackrel{\mathrm{d}}{=}}

\numberwithin{equation}{section}

\endlocaldefs

\begin{document}

\begin{frontmatter}
\title{PageRank Asymptotics on Directed Preferential Attachment Networks}
\runtitle{PageRank}

\begin{aug}
\author{\fnms{Sayan} \snm{Banerjee}\ead[label=e1]{sayan@email.unc.edu}}
\and
\author{\fnms{Mariana} \snm{Olvera-Cravioto}\ead[label=e2,mark]{molvera@email.unc.edu}}

\address{Statistics and Operations Research, UNC-Chapel Hill, \printead{e1,e2}}
\end{aug}

\begin{abstract}
We characterize the tail behavior of the distribution of the PageRank of a uniformly chosen vertex in a directed preferential attachment graph and show that it decays as a power law with an explicit exponent that is described in terms of the model parameters. Interestingly, this power law is heavier than the tail of the limiting in-degree distribution, which goes against the commonly accepted {\em power law hypothesis}. This deviation from the power law hypothesis points at the structural differences between the inbound neighborhoods of typical vertices in a preferential attachment graph versus those in static random graph models where the power law hypothesis has been proven to hold (e.g., directed configuration models and inhomogeneous random digraphs). In addition to characterizing the PageRank distribution of a typical vertex, we also characterize the explicit growth rate of the PageRank of the oldest vertex as the network size grows.
\end{abstract}



\end{frontmatter}
\section{Introduction}

PageRank, the ranking scheme introduced by Brin and Page \cite{page1999pagerank}, is arguably one of the most important centrality measures for directed complex networks. Originally intended to assign a universal rank to pages in the World Wide Web, it has become a standard tool in the analysis of almost any large graph, both directed or undirected. PageRank's appeal is due in part to its computability, which involves solving a linear system of equations that can be efficiently done in a distributed fashion even on very large graphs \cite{sarma2015fast}. But more importantly, it is PageRank's ability to identify ``influential" nodes that has made it popular in a wide range of applications \cite{agirre2014random, Andersen06,   Chen06citations,  gleich2015pagerank, Gyongyi04, Haveliwala2002personalization, ivan2011web,jing2008visualrank, waltman2010eigenfactor, wang2007keyword}. It is this vague concept of ``influence" that has been at the core of the analysis of PageRank's qualitative behavior, whose aim is to understand the type of nodes that PageRank will tend to score highly. 

Since the early 2000's, research on the qualitative analysis of PageRank has focused on characterizing the distribution of the ranks it produces in relation to the properties of the underlying graph. The work in \cite{pandurangan2002using} and other papers that followed \cite{becchetti2006distribution, donato2004large, fortunato2006approximating}, identified a property known as the {\em power law hypothesis}, which states that in a directed network whose in-degree distribution follows a power law\footnote{$X$ has a power law distribution if $P(X > x)$ decays proportionally to $x^{-\alpha}$ for some $\alpha > 0$.}, the PageRank scores will also follow a power law with the same exponent. This phenomenon, that has been empirically verified in many real-world networks, has been proven to be true for directed random graphs that are generated via either a directed configuration model \cite{chen2017generalized, olvera2019pagerank} or an inhomogeneous random digraph \cite{lee2020pagerank, olvera2019pagerank}, and is based on the asymptotic analysis of the local weak limit of the underlying graphs \cite{jelenkovic2010information, olvera2012tail, volkovich2010asymptotic}. Both of these random graph models are {\em static}, i.e., they are meant to describe a fixed instance of a real-world graph. The two models share other important characteristics, such as their typically short distances and diameter (known as the {\em small world phenomenon}), their ability to replicate almost any given degree distribution, in particular power laws (known as the {\em scale-free} property), and their local tree-like structure \cite{van2016random}. In addition to establishing the validity of the power law hypothesis for the two models, the analysis of the PageRank distribution done in \cite{jelenkovic2010information, olvera2019pagerank} provides insights into the types of nodes that PageRank will tend to rank highly, which broadly classifies highly ranked nodes into two categories: 1) nodes that have an atypically large in-degree, and 2) nodes that have one inbound neighbor with an unusually large PageRank.  See also the recent work \cite{cai2021rankings} for a large deviations type comparison between the maximum in-degree and the maximum PageRank on the directed configuration model.

This paper further contributes to the qualitative analysis of PageRank, however, on a different kind of random graphs. A popular alternative for modeling real-world networks is to consider {\em evolving} random graph, which not only mimic static properties such as the small-world phenomenon or the scale-free property, but also provide an explanation for how the graphs are created. Perhaps the most influential of such models is the preferential attachment graph proposed by Albert and Barab\'asi \cite{barabasi1999emergence}. Here we focus on a directed preferential attachment model (DPA), in which vertices are added to the graph one at a time by drawing directed edges towards $m$ vertices chosen according to a linear preference rule.  This model is known to produce graphs with power law in-degrees, small diameter, small typical distances, and a local tree-like structure. The main difference compared to its static counterparts is that it also provides a time-stamp identifying when each vertex was added to the graph. Interestingly, as will be shown in the current article, this time-stamp has important consequences for the qualitative behavior of PageRank, which illustrates how different the neighborhoods of large in-degree vertices are relative to those in the static models.  More precisely, our work shows that in the DPA, the PageRank distribution follows a power law, however, the index of this power law is different from that of the in-degree distribution, with the former having even heavier tails. 

The expected value of PageRank on DPA graphs was first studied in \cite{avrachenkov2006pagerank}, for the same model studied here. In addition to the expected values, the work in \cite{avrachenkov2006pagerank} also hints at the power law distribution of PageRank, although stops short of proving it. More recently, the work in \cite{garavaglia2020local} uses a general approach based on local weak convergence to establish a lower bound for the tail distribution of PageRank in terms of that of the in-degree, which combined with the known local weak convergence of the DPA gives a power law lower bound. Our first main result, which is enabled by a novel reinterpretation of the local weak limit of the DPA, shows that the limiting PageRank distribution of a typical vertex is indeed a power law with an explicit index that is different from that of the in-degree distribution. In other words, the DPA provides a counterexample to the power law hypothesis. Our second main result furnishes precise asymptotics of the PageRank of the oldest vertex in the network. The proofs of the main results involve a careful analysis of the local weak limit of the DPA, which we show can be expressed in terms of a continuous time branching process. This branching process can then be used to construct a continuous time martingale from where the properties of the PageRank distribution can be derived. 

The DPA model considered here is limited in the sense that it assigns to each vertex a deterministic out-degree and it produces graphs having no directed cycles, as edges always point from younger to older vertices. A variant of our model capable of producing directed cycles and random out-degrees is the one proposed in \cite{bollobas2003directed}, in which at each time step,  either a new vertex attaches to the existing network using an edge with a random direction, or a new directed edge is added between existing vertices. The asymptotic joint in-degree and out-degree distribution of a uniformly chosen vertex in this model was recently analyzed in \cite{britton2018directed}, in a unified setting that also includes more realistic versions of the model. Since we believe that our techniques can be broadly applied to other random network models whose local weak limit can be described in terms of continuous time branching processes, we expect to extend our results to other DPA models in the future.

This paper is organized as follows. In Section~\ref{modeldef} we define the DPA and give a mathematical description of PageRank. In Section~\ref{main} we state our main results on the asymptotic behavior of PageRank and compare it to its counterparts on static random graph models. We also discuss the modeling implications of our results, which provide important insights about the local neighborhoods of high degree vertices. Section~\ref{locsec} summarizes known results about the local weak convergence of the DPA. Section~\ref{loclimunif} contains a new theorem unifying two different representations of the local limit, a key step for the analysis of PageRank.  Finally, Sections~\ref{proofs} and \ref{rootsec} contain the proofs of our main results.

\section{PageRank and the directed preferential attachment model}\label{modeldef}
PageRank is a centrality measure which gives to each vertex in a graph a score based on the stationary distribution of a discrete-time Markov chain on the graph, as described below.

Let $G = G(V, E)$ denote a directed graph with set of vertices $V $ and set of edges $E$. Since the graph $G$ is allowed to have self-loops and multiple edges from one vertex to another, we use $a_{u,v}$ to denote the number of edges in $G$ from vertex $u$ to vertex $v$. For each vertex $v \in V$ let $d_v^-$ and $d_v^+$ denote its in-degree and out-degree, respectively. Writing $|V|$ for the number of vertices in the graph, let $A$ denote the adjacency matrix of $G$ defined as the $|V| \times |V|$ matrix whose $(i,j)$th element is $a_{i,j}$, and let $\Delta$ be the diagonal matrix whose $i$th element is $1/d_i^+$ if $d_i^+>0$ and $0$ if $d_i^+=0$. Let $P$ be the matrix product $\Delta A$ with the zero rows replaced with the probability vector ${\bf q} = |V|^{-1} \mathbf{1}$. Note that $P$ is a stochastic matrix (all rows sum to one). The PageRank vector $\boldsymbol{\pi} = (\pi_1, \dots, \pi_{|V|})$, with \emph{damping factor} $c \in (0,1)$, is defined as the stationary distribution of the discrete-time Markov chain which, at each time, hops according to transition matrix $P$ with probability $c$, and jumps to a uniformly chosen vertex in $G$ with probability $1-c$. In particular, when this Markov chain is at a vertex $v$ with $d_v^+>0$, the hop mechanism lands it on a neighboring vertex via an incident edge chosen uniformly at random. When the chain is at a \emph{dangling node} (with out-degree $0$), its next move is always to a uniformly chosen vertex.

The PageRank vector $\boldsymbol{\pi}$ can be computed by solving the following system of equations \cite[Section 1.5]{levin2017markov}:
$$\boldsymbol{\pi} = \boldsymbol{\pi} (cP) + (1-c) \mathbf{q}.$$
As the matrix $cP$ is substochastic (its rows sum to $c$), the system of equations is guaranteed to have a unique solution given by:
$$\boldsymbol{\pi} = (1-c)\mathbf{q}(I - cP)^{-1} = (1-c) \mathbf{q} \sum_{k=0}^\infty (c P)^k.$$

For the purpose of analyzing the typical behavior of PageRank, it is more convenient to work with the {\em scale-free} PageRank vector $\mathbf{R} := |V| \boldsymbol{\pi}$. The object of our study is a uniformly chosen component of $\mathbf{R}$. 

\bigskip

We now move on to describe the random graph model that is the focus of this paper. The directed (linear) preferential attachment model $\DPA(m,\beta)$ with parameters $m \in \mathbb{N}$ and $\beta \ge 0$ is a model for an evolving random directed rooted graph sequence $\{\G_n\}_{n \ge 0}$ obtained through the following recipe. $\G_1$ comprises one vertex $v_1$ (the root) and zero edges. $\G_2$ contains two vertices $v_1$ and $v_2$ connected by $m$ edges directed from $v_2$ to $v_1$. Given we have obtained $\G_{n-1}$ for $n \ge 3$, $\G_n$ is constructed from $\G_{n-1}$ by adding one vertex $v_n$ with $m$ outbound edges which are connected one-by-one to the existing vertices $\{v_i\}_{1 \le i \le n-1}$ with probability proportional to their degree, i.e., for $1 \le k \le m, 1 \le i \le n-1$:
$$
\mathbb{P}\left( \left. k^{th} \text{ outbound edge of } v_n \text{ is attached to } v_i \right| \G_{n-1}\right) = \frac{D_i(n-1,k-1) + \beta}{\sum_{j=1}^{n-1}(D_j(n-1,k-1) + \beta)},
$$
where $D_i(n-1,k-1)$ is the total degree (in-degree plus out-degree) of the $i^{th}$ vertex after $k-1$ edges of $v_n$ are attached to vertices in $\G_{n-1}$. Note that the denominator simplifies to $\sum_{j=1}^{n-1}(D_j(n-1,k-1) + \beta) = 2m(n-2) + k-1 + \beta(n-1)$. 

Note that the in-degree of vertices in the graph sequence may grow with time, while their out-degree is always equal to $m$. Denoting the adjacency matrix of $\G_n$ by $A_n$, the constant out-degree also implies that the matrix $cP$ simplifies to $(c/m) A_n$, and the scale-free PageRank vector $\mathbf{R}(n)$ on $\G_n$ becomes:
$${\bf R}(n) = (c/m) \mathbf{ R}(n)  A_n + (1-c) {\bf 1} = (1-c)  \sum_{k=0}^\infty  \left( \frac{c}{m} \right)^k {\bf 1} A_n^k.$$
It will also be convenient in the sequel to introduce the notation:
$$\mathcal{P}_{k,i}^{(n)} = (\mathbf{1} A_n^k)_i \qquad k \geq 1,$$
which corresponds to the number of directed paths of length $k$ in $\G_n$ that end at vertex $i$. From now on, we will also denote the PageRank vector by $(R_1(n), \dots, R_n(n))$. 
Hence, 
we obtain the following representation of PageRank for the $\DPA(m,\beta)$:
\begin{equation}\label{pagedef}
R_i(n) = (1-c)\left(1 + \sum_{k=1}^{\infty}\left(\frac{c}{m}\right)^k \PP^{(n)}_{k,i}\right), \quad i = 1, \dots, n.
\end{equation}

\section{Main Results}\label{main}

Our main result, and the one that relates to the power law hypothesis, answers a question about the behavior of a {\em typical} vertex,  represented by a uniformly chosen vertex in $\G_n$. However, since the DPA is an evolving graph model where vertex labels tell us their age, another interesting quantity to analyze is the behavior of the oldest vertex, which is quite different from that of a typical vertex. Hence, this section is divided into two subsections, one for each of these two cases.

\subsection{PageRank of a uniformly chosen vertex} 

Let $V_n$ denote (the index of) a uniformly chosen vertex in $\G_n$. Denote the in-degree and  PageRank of this vertex by $D^{-}_{V_n}(n)$ and $R_{V_n}(n)$, respectively. 
The following theorem shows that these quantities jointly converge in distribution as $n \rightarrow \infty$. Moreover, it explicitly quantifies the power law tail exponents of the limiting in-degree and that of the PageRank distribution. 
\begin{theorem}\label{exponent}
The in-degree and PageRank of a uniformly chosen vertex in $\G_n$ jointly converge in distribution:
\begin{equation}\label{weakcon0}
(D_{V_n}^-(n), R_{V_n}(n)) \xrightarrow{d} (\mathcal{D}^-, \mathcal{R}) \ \text{ as } n \rightarrow \infty,
\end{equation}
where $\mathcal{R}, \mathcal{D}^-$ are described in Section \ref{limdr}. Moreover, there exists $C>0$ such that as $k \rightarrow \infty$,
\begin{equation}\label{degtail}
\pr{\mathcal{D}^- \ge k} = Ck^{-2-\beta/m}\left(1 + O(k^{-1})\right),
\end{equation}
and positive constants $C_1, C_2$ such that for any $r \geq 1$, 
\begin{equation}\label{prtail}
C_1r^{-(2+\beta/m)/(1+(m+\beta)c/m)} \ \le \ \pr{\mathcal{R} \ge r} \ \le \ C_2r^{-(2+\beta/m)/(1+(m+\beta)c/m)}.
\end{equation}
\end{theorem}

\begin{remark}\label{powerdiff}
Theorem \ref{exponent} leads to the following observations.

(i) The limiting PageRank distribution has a \emph{heavier tail} than the limiting in-degree distribution. This stands in stark contrast with the commonly accepted power law hypothesis which asserts that, in graphs whose in-degree follows a power law, the tail distribution of PageRank also follows a power law with the same exponent. A more thorough discussion on the implications of this result is included at the end of the subsection.

(ii) The asymptoic behavior of the limiting in-degree given by \eqref{degtail} shows that its tail exponent is linearly increasing in $\beta$ (keeping $m,c$ fixed) and thus the tails become lighter as $\beta$ grows. This can be intuitively understood by noting that increasing $\beta$ makes the attachment mechanism `approach' uniform attachment (each new vertex attaches to a pre-existing vertex uniformly at random) where the limiting degree distribution has exponential tails \cite[Theorem 1.1]{janson2005asymptotic}. 

However, the tail exponent of the limiting PageRank distribution in \eqref{prtail} satisfies:  
$$
\lim_{\beta \rightarrow \infty} \frac{2+\beta/m}{1+(m+\beta)c/m} = \frac{1}{c},
$$
which shows that PageRank remains a power law as $\beta$ increases, with a tail exponent uniformly bounded above and below by positive numbers. In fact, this suggests the surprising phenomenon that the limiting PageRank distribution in the uniform attachment model has a power law tail with exponent $1/c$, although the limiting in-degree distribution has exponential tails. Although we do not prove it in this paper, we believe this is true. In the tree case $(m= 1)$, this can be shown directly using the methods developed in the current paper. The non-tree case requires a description of the local weak limit in the uniform attachment setting, and will be included in future work involving a more general class of dynamic random graphs.

(iii) As $c \downarrow 0$ (keeping $m,\beta$ fixed), the tail exponent of the limiting PageRank distribution approaches that of the limiting in-degree distribution. This can be intuitively understood from the representation of PageRank in \eqref{pagedef} by noting that as $c \downarrow 0$ the main contribution to the PageRank of any vertex comes from its number of inbound neighbors, which is precisely the in-degree of this vertex.

(iv) As $m \rightarrow \infty$ (keeping $\beta,c$ fixed), the tail exponent of the limiting in-degree distribution approaches $2$ (which also corresponds to the $\beta=0$ case). However, the exponent for the limiting PageRank distribution approaches $2/(1+c)$.
\end{remark}

{\em The power law hypothesis.} Although numerous studies  \cite{becchetti2006distribution, donato2004large, fortunato2006approximating, pandurangan2002using} have empirically verified the validity of the power law hypothesis on real-world scale-free graphs, Theorem~\ref{exponent} provides a counterexample, since although it establishes the power law behavior of both the limiting PageRank distribution and that of the limiting in-degree, the two tail exponents are different. 

To shed some light into the implications of our results, it may be helpful to explain in more detail what the existing theorems for static random graphs say. First, both the directed configuration model and the inhomogeneous random digraph, converge locally in the large graph limit to a (possibly infinite) marked Galton-Watson process. This implies that the limiting PageRank distribution can be characterized through a branching distributional fixed-point equation. In the case of the directed configuration model with in-degree distributed according to $\mathcal{D}^-$ in Theorem~\ref{weakcon0} and out-degree equal to $m$ for all nodes, Theorem 6.4 in \cite{chen2017generalized} characterizes the limiting PageRank distribution as the endogenous solution to:
\begin{equation} \label{fixedpoint} 
\mathcal{R} \stackrel{d}{=} \sum_{i=1}^{\mathcal{D}^-} \frac{c}{m} \mathcal{R}_i + 1-c,
\end{equation}
where the $\{ \mathcal{R}_i\}$ are i.i.d.~copies of $\mathcal{R}$, independent of $\mathcal{D}^-$. In other words, the PageRanks of the inbound neighbors of the randomly chosen vertex are asymptotically i.i.d.~and independent of its in-degree, which leads to the heavy-tailed asymptotic 
\begin{align*}
\pr{\mathcal{R} > x} & \sim \pr{\mathcal{D}^- > (c \mathbb{E}[\mathcal{R}]/m)^{-1} x } + \pr{ \max_{1 \leq i \leq \mathcal{D}^-} \mathcal{R}_i > (m/c) x } \\
&\sim C \pr{ \mathcal{D}^- > x }, \qquad  x  \to \infty.
\end{align*}
The second equivalence is \cite[Theorem 5.1]{jelenkovic2010information} while the first is a consequence of the proof of the same theorem\footnote{$g(x) \sim f(x)$ as $x \rightarrow \infty$ if $\lim_{x \rightarrow \infty} f(x)/g(x) =1$.}.
On the other hand, the DPA converges in the local weak sense to a continuous-time branching process stopped at a finite random time, that is, a finite tree. Moreover, the time-stamps imply that the PageRanks of the inbound neighbors of a randomly chosen vertex are no longer i.i.d., nor are they independent of its in-degree. In other words, the distributional fixed-point equation \eqref{fixedpoint} does not hold, and the asymptotic behavior of the random variable $\sum_{i=1}^{\mathcal{D}^-} (c/m) \mathcal{R}_i + 1-c$ is more complex (see also \cite[Remarks 6.7 and 6.11]{garavaglia2020local}). The heavier tails of $\mathcal{R}$, relative to those of $\mathcal{D}^-$, are consistent with the observation that ``old" vertices in the DPA (strongly correlated with having large in-degrees) tend to have inbound neighbors that are also ``old", hence compounding the effect of their large in-degrees.

In other words, the tail behavior of PageRank on the DPA, as stated by Theorem~\ref{weakcon0}, reflects the different structure of local neighborhoods in the DPA compared to local neighborhoods in either a directed configuration model or an inhomogeneous random digraph. In the former, large degree vertices tend to be close to each other, while in the latter, they are more evenly spread out. Moreover, PageRank can be used to statistically distinguish the DPA from either of the two static models for which the power law hypothesis has been shown to hold.

\subsection{PageRank asymptotics for the oldest vertex} The next theorem quantifies the asymptotic behavior of the PageRank of the oldest vertex in the $\DPA(1,\beta)$, which corresponds to the root of a branching process (in the $m=1$ case, the network $\G_n$ is a tree for all $n \in \mathbb{N}$). As our theorem shows, the behavior of the oldest vertex is quite different from that of a typical vertex.

\begin{theorem}\label{rootpr}
The PageRank $R_1(n)$ of the root in $\G_n$ for the $\DPA(1,\beta)$ model satisfies
\begin{equation}\label{rootas}
n^{-(1+(1+\beta)c)/(2+\beta)} R_1(n) \xrightarrow{a.s} W, \qquad  n\to \infty,
\end{equation}
for some positive, almost surely finite random variable $W$.
\end{theorem}

\begin{remark}\label{persist}
Theorem \ref{rootpr} leads to the following natural questions. Although we do not explore these in the current paper, we will address them in future research.
\begin{itemize}
\item[(i)] \emph{Root PageRank asymptotics for $m \ge 2$.} The proof of Theorem \ref{rootpr} relies on an embedding of the discrete tree network in a continuous time branching process (see Lemma \ref{lem:ctb-embedding}). Although such an embedding does not directly extend to the non-tree case, there is a natural way of obtaining $\DPA(m,\beta)$ from $\DPA(1,\beta)$ by \emph{collapsing a tree network} (equivalently, the associated continuous time branching process) \cite{garavaglia2018trees}. Broadly speaking, this procedure gives a $\DPA(m,\beta)$ network of size $n$ from a $\DPA(1,\beta)$ tree of size $nm$ by collapsing $m$ successive vertices, in their arrival order, along with their incident edges. However, vertices at different distances from the root in the tree can be collapsed into the same vertex in the $\DPA(m,\beta)$ network, which makes the analysis of PageRank hard to transfer from the tree to the non-tree case. Hence, this case remains an open problem.
\item[(ii)] \emph{Maximal PageRank.} We believe the asymptotic growth rate of the root PageRank is of the same order as that of the \emph{maximal PageRank} of the network. This belief stems from the phenomenon of \emph{persistence} of the maximal degree vertex: almost surely, there exists a random $n_0$ such that the maximal degree vertex in $\G_{n_0}$ continues to be the maximal degree vertex in $\G_n$ for all $n \ge n_0$ \cite{banerjee2020persistence, galashin2013existence}. In other words, we conjecture that the maximal PageRank also exhibits persistence, which would imply that Theorem~\ref{rootpr} gives its asymptotic behavior. 
\end{itemize}
\end{remark}

\section{PageRank and Local Weak Convergence}\label{locsec}

The recent work in \cite{garavaglia2020local} shows that the PageRank $R_{V_n}(n)$ of a uniformly chosen vertex in $\G_n$ converges in distribution to a random limit $\mathcal{R}$ that is explicitly characterized by the local weak limit of $\G_n$ as $n \rightarrow \infty$. We will now describe this connection, which will establish \eqref{weakcon0}.

\subsection{Local weak convergence for directed graphs}

Local weak limits characterize the asymptotic behavior of neighborhoods of a uniformly chosen vertex in a graph sequence as the size of the graph goes to infinity. This concept was first introduced in \cite{aldous2007processes,aldous2004objective,benjamini2011recurrence} for undirected graphs. We now sketch an extension of this concept to directed graphs as laid out in \cite{garavaglia2020local}. 

Let $\mathbb{G}$ denote the space of directed, marked, rooted graphs \cite[Definition~3.8]{garavaglia2020local}. Elements of this space comprise directed graphs $G$ with a distinguished vertex $\emptyset$ called the root. Moreover, each vertex carries an integer value called the mark which is at most the out-degree of the vertex\footnote{ For the $\DPA$ models considered in this paper, all marks can be taken to be $m$ as each vertex (other than the root) has the same out-degree. However, we will work in the marked setting to apply some results of \cite{garavaglia2020local} without modification.}. There is a natural concept of isomorphism `$\cong$' between two elements of $\mathbb{G}$ \cite[Definition 3.9]{garavaglia2020local}: two such elements are isomorphic if there exists a bijection between the vertex sets which maps root to root and preserves the directed adjacency structure and marks of the vertices. We will denote by $\mathbb{G}_{\star}$ the quotient space of $\mathbb{G}$ with respect to the equivalence given by isomorphisms. We will denote a generic element of $\mathbb{G}_{\star}$ by $(G, \emptyset, M(G))$, where $M(G)$ denotes the set of marks on vertices of $G$.

\textbf{Notation:} We will label the vertices of trees and individuals in branching processes by $\II := \cup_{k=0}^{\infty} \mathbb{N}^k$, with the convention $\mathbb{N}^0 := \{\emptyset \}$. In addition, for $\oa = (i_1, \dots, i_k)$, we denote by $(\oa, j) = (i_1, \dots, i_k, j)$ the index concatenation operation, and $|\oa| = k$ the length of the label, or equivalently, the generation to which the individual in the branching process belongs to. For simplicity, for labels of length one (equiv. individuals in the first generation), we omit the parenthesis and simply write $i \in \mathbb{N}$. Furthermore, for a vector $\oa$ of length $|\oa| = k$, we use the notation $\oa|m = (i_1, \dots, i_m)$ for $0 \leq m \leq k$ to denote its truncation to length $m$, with the convention $\oa|0 = \emptyset$. For continuous time branching processes, the birth time of an individual $\oa \in \II$ will be denoted by $\sigma_\oa$. For $i \in \mathbb{N}$, we will denote by $\mathbf{e}^{(i)}$ the $i$-th coordinate unit vector in $\mathbb{R}^{\infty}$.  For $n \in \mathbb{N}$ and any $n \times n$ matrix \textcolor{black}{$Q$, $e^{Q}$} will denote the usual matrix exponential.

For any $k \in \mathbb{N} \cup \{0\}$, the $k$-neighborhood of the root $\emptyset$ in $(G, \emptyset, M(G))$, denoted by $U_{\le k}(\emptyset)$, is obtained by progressively exploring vertices using incoming edges \emph{in the opposite direction}, starting from the root, up till graph distance $k$ from the root and revealing the marks and connectivity structure of the explored vertices \cite[Definition 3.10]{garavaglia2020local}. Note that $U_{\le k}(\emptyset)$ thus constructed is a marked directed subgraph of $(G, \emptyset, M(G))$. This leads to a natural distance on the space $\mathbb{G}_{\star}$: for two elements $(G, \emptyset, M(G))$, $(G', \emptyset', M(G'))$ in $\mathbb{G}_{\star}$, define $d_{loc}((G, \emptyset, M(G)), (G', \emptyset', M(G'))) := 1/(\kappa + 1)$, where $\kappa := \inf\{k \ge 1 : U_{\le k}(\emptyset) \not\cong U_{\le k}(\emptyset')\}$. Unlike for undirected graphs, this distance is not a metric on $\mathbb{G}_{\star}$, but a pseudonorm, as the distance between two distinct elements in $\mathbb{G}_{\star}$ can be zero. This is because, in the directed setting, edges are explored only in one direction, leaving parts of the graph unexplored \cite[Figure 4]{garavaglia2020local}. Thus, two rooted marked directed graphs can have identical explorable root neighborhoods without being isomorphic. However, if $d_{loc}((G, \emptyset, M(G)), (G', \emptyset', M(G')))=0$, it can be shown that the incoming neighborhoods for the root (the possibly infinite subgraph that can be explored starting from the root) in the two graphs are isomorphic. Namely, denoting the respective incoming neighborhoods by $U_{\infty}(\emptyset)$ and $U_{\infty}(\emptyset')$, we have $U_{\infty}(\emptyset) \cong U_{\infty}(\emptyset')$. This leads to an equivalence relation on $\mathbb{G}_{\star}$. Denoting by $\tilde{\mathbb{G}}_{\star}$ the quotient space of $\mathbb{G}_{\star}$ under this equivalence relation, it follows that $(\tilde{\mathbb{G}}_{\star}, d_{loc})$ is a Polish space.

Now, we have all the tools to describe the concept of local weak convergence in the directed graph setting. For a sequence $\{(G_n, M(G_n))\}_{n \in \mathbb{N}}$ of marked, directed random graphs, define 
$$
\mathbb{P}_n := \frac{1}{|V(G_n)|} \sum_{v \in V(G_n)} \delta_{(G_n, v, M(G_n))}
$$
as the empirical measure corresponding to selecting the root in $(G_n, M(G_n))$ uniformly at random in $V(G_n)$.

\begin{definition}[Local weak convergence for directed graphs]
Consider a sequence $\{(G_n, M(G_n))\}_{n \in \mathbb{N}}$ of marked, directed random graphs. We say $(G_n, M(G_n))$ converges \emph{in probability in the local weak sense} to a random element $\mathcal{G}^* \in \tilde{\mathbb{G}}_{\star}$ with law $\mathbb{P}^*$ if for any bounded continuous function $f : \tilde{\mathbb{G}}_{\star} \rightarrow \mathbb{R}$,
$$
\mathbb{E}_{\mathbb{P}_n}(f) \xrightarrow{P} \mathbb{E}_{\mathbb{P}^*}(f) \qquad n \to \infty,
$$
where $\mathbb{E}_{\mathbb{P}_n}$ and $\mathbb{E}_{\mathbb{P}^*}$ respectively denote expectation taken with respect to the laws $\mathbb{P}_n$ and $\mathbb{P}^*$.
\end{definition}

\subsection{Local weak limits of $\DPA(m,\beta)$}\label{loclimsec}

We will now describe local weak limits of directed linear preferential attachment graphs. The limiting graphs in the cases $m=1$ and $m \ge 2$ have both been described before in \cite{rudas-2} and \cite{berger2014asymptotic}, respectively, but the descriptions are very different. A key technical contribution of this article is a unified description in terms of continuous time branching processes that is valid for $m\geq 1$. In this subsection, we offer the existing descriptions of the limiting graphs for the two cases $m=1$ and $m \ge 2$. The
new unified description is given in Section \ref{loclimunif}.

(I) \textit{Local limit for $m=1$ case}. A few definitions are in order.

\begin{definition}[$\beta$-Yule Process]\label{Yule}
Fix $\beta \ge 0$. A $\beta$-Yule process is a pure birth process $\{Y_{\beta}(t):t\geq 0\}$ with $Y_{\beta}(0)=0$ and which satisfies, for any $t \ge 0$, 
\begin{align*}
\pr{Y_{\beta}(t+dt) - Y_{\beta}(t) = 1|\mathcal{F}(t)} &:= (Y_{\beta}(t) + 1 + \beta) dt + o(dt) \qquad \text{and} \\
\pr{Y_{\beta}(t+dt) - Y_{\beta}(t) \ge 2|\mathcal{F}(t)} &:= o(dt),
\end{align*}
where $\{\mathcal{F}(t):t\geq 0\}$ is the natural filtration of the process. 
\end{definition}

\begin{definition}[Continuous time branching process (CTBP)]\label{ctbp}
Fix $\sigma_\emptyset \ge 0$. Let $\{\xi_x(t) : t \in [x, \infty)\}, x \in [\sigma_\emptyset, \infty),$ be a collection of independent point processes indexed by the real numbers in $[\sigma_\emptyset,\infty)$. A continuous time branching process (CTBP) with root $\emptyset$ born at time $\sigma_\emptyset$, denoted by $\{\B_{\sigma_\emptyset}(t) : t \ge \sigma_\emptyset \}$, is a branching process which originates from a single individual $\emptyset \in \mathcal{I}$ which gives birth to new individuals indexed by $\{i : i \in \mathbb{N}\}$ at times prescribed by the point process $\xi_{\sigma_\emptyset}(\cdot)$. Moreover, any individual $\oa \in \II$ born into the population at time $\sigma_\oa$ produces offspring $\{(\oa,j) : j \in \mathbb{N}\}$ independently at times given by $\xi_{\sigma_\oa}(\cdot)$. The associated directed rooted tree $\TT(\B_{\sigma_\emptyset}(t))$ for the CTBP observed till time $t$ is constructed by placing an edge between each individual and its parent which is directed towards the parent. The total number of individuals in the CTBP at time $t$ will be denoted by $|\B_{\sigma_\emptyset}(t)|$. 
\end{definition}

We now describe the local weak limit for $\DPA(1,\beta)$ for fixed $\beta \ge 0$ (all vertex marks are equal to 1). In this case, the local weak limit of $\G_n$ as $n \rightarrow \infty$ is described in terms of a CTBP with root $\emptyset$ born at time $0$ where each individual $\oa \in \II$ reproduces according to the point process $\xi_{\sigma_\oa}(t) = Y_\oa(t-\sigma_\oa)$, $t \ge \sigma_\oa$, where $\{Y_\oa(\cdot) : \oa \in \II\}$ is an i.i.d.~collection of $\beta$-Yule processes. We denote this CTBP by $\{\B^{(\beta)}(t): t \ge 0\}$. It was shown in \cite{rudas-2} (see Theorem 1 and Section 4.2 there) that the local weak limit of $\G_n$ as $n \rightarrow \infty$ is given by $\G^*(1) := \TT\left(\B^{(\beta)}(\tau)\right)$, where $\tau$ is an $\operatorname{Exp}(2+\beta)$ random variable that is independent of $\B^{(\beta)}(\cdot)$.

(II) \textit{Local limit for $m \ge 2$ case}. We now describe a random tree, called the \emph{P\'olya point graph}, that arises as the local weak limit of $\DPA(m,\beta)$ in the $m \ge 2$ case (all vertex marks are equal to $m$). This graph was defined in \cite{berger2014asymptotic} (see also \cite[Definition 6.9]{garavaglia2020local}).

\begin{definition}[P\'olya point graph]\label{polya}
Fix $m \ge 2$. Let $\{\gamma_{x} : x \in (0, \infty)\}$ be a collection of i.i.d.~$\operatorname{Gamma}(m+\beta,1)$ random variables indexed by the positive real numbers. Set $\chi := (m+\beta)/(2m+\beta)$ and $\psi:= (1-\chi)/\chi$. The P\'olya point graph $\G^*(m)$ is defined as the rooted random tree defined as follows:
\begin{itemize}
\item[(i)] The root $\emptyset$ is assigned a birth time $\sigma_\emptyset \sim U^{\chi}$, where $U \sim \operatorname{Uniform}(0,1)$. The root reproduces according to the \textcolor{black}{nonhomogeneous} Poisson point process with intensity
$$
\rho_{\sigma_\emptyset}(v) := \gamma_{\sigma_\emptyset}\frac{\psi v^{\psi-1}}{\sigma_\emptyset^{\psi}}, \ v \in [\sigma_\emptyset,1].
$$
Denote the offsprings by $\{ i : 1 \le i \le N_\emptyset\}$, where $N_\emptyset$ is the total number of births in $[\sigma_\emptyset,1]$. Each offspring is attached to the root by an edge directed towards the root.
\item[(ii)] Each individual $\oa \in \II$ with $\sigma_\oa \le 1$ reproduces independently according to the  \textcolor{black}{nonhomogeneous} Poisson point process 
$$
\rho_{\sigma_\oa}(v) := \gamma_{\sigma_\oa}\frac{\psi v^{\psi-1}}{\sigma_\oa^{\psi}}, \ v \in [\sigma_\oa,1].
$$
Its offsprings are denoted by $\{(\oa,j) : 1 \le j \le N_\oa\}$, where $N_\oa$ is the total number of births in $[\sigma_\oa, 1]$. A directed edge pointing towards $\oa$ is placed between $\oa$ and each of its offsprings.
\end{itemize}
\end{definition}

\begin{remark}\label{polyactbp}
Note that the P\'olya point graph $\G^*(m)$ is obtained as $\TT\left(\B^{pol}_{\sigma_\emptyset}(1)\right)$, where $\sigma_\emptyset \sim U^{\chi}$ for a $\operatorname{Uniform}(0,1)$ random variable $U$ and, conditionally on $\sigma_\emptyset = z$, $\{\B^{pol}_{z}(t) : t \ge z\}$ is a CTBP with associated reproduction point process for each $x \in [z, \infty)$ given by a  \textcolor{black}{nonhomogeneous} Poisson point process $\{\xi^{pol}_x(t) : t \in [x, \infty)\}$ with (random) intensity $\rho_x(v) := \gamma_{x}\frac{\psi v^{\psi-1}}{x^{\psi}}, \ v \in [x,\infty)$.
\end{remark}

\subsection{Limiting joint distribution of PageRank and in-degree of a uniformly chosen vertex} \label{limdr}
The representation \eqref{pagedef} shows that the dependence of the PageRank of a vertex on the graph structure outside a neighborhood of the vertex `diminishes exponentially' with the size of the neighborhood. From this, it is plausible that the distributional limit of the PageRank $R_{V_n}(n)$ of a uniformly chosen vertex indexed $V_n$ in $\G_n$ can be described in terms of statistics on the graph that appears as the local weak limit of $\G_n$ as $n \rightarrow \infty$. Indeed, such a heuristic  was formally justified in \cite{garavaglia2020local} for any collection of random graphs that converge in the local weak sense \cite[Theorem~2.1]{garavaglia2020local}. We take this approach for the directed preferential attachment model and slightly extend the result of \cite{garavaglia2020local} to describe the joint limiting in-degree and PageRank of a uniformly chosen vertex.

It is easy to see that $\G^*(m)$ is almost surely finite for each $m \ge 1$. Define the random variable
\begin{equation}\label{limpr}
\mathcal{R} := (1-c)\left(1 + \sum_{k=1}^{\infty}\left(\frac{c}{m}\right)^k \PP^*_k\right)
\end{equation}
where $\PP^*_k$ denotes the number of directed paths of length $k$ in $\G^*(m)$ that end at the root. Also, denote by $\mathcal{D}^-$ the in-degree of the root in $\G^*(m)$. For notational convenience, we suppress the dependence of $\mathcal{R},\mathcal{D}^-$ on $m$.

The following convergence assertion is a straightforward extension of \cite[Theorem 2.1]{garavaglia2020local}, however, we include its proof to add some steps that are missing in \cite{garavaglia2020local}. Recall that $D^-_{V_n}(n)$ and $R_{V_n}(n)$ denote, respectively, the in-degree and the PageRank of a uniformly chosen vertex indexed $V_n$ in $\G_n$. 
\begin{theorem}\label{weakconv}
Fix $m \in \mathbb{N}$. Then for any continuity point $r$ of the distribution function of $\mathcal{R}$ and $k \in \mathbb{N} \cup \{0\}$, we have
\begin{equation}\label{prcon}
\frac{1}{n}\sum_{i=1}^n\mathbb{I}(D^-_i(n) \ge k, R_i(n) > r) \xrightarrow{P} \mathbb{P}(\mathcal{D}^- \ge k, \mathcal{R} > r)
\end{equation}
as $n \to \infty$. In particular, \eqref{weakcon0} holds.
\end{theorem}

\begin{proof}
Let $\mathcal{G}_n$ denote a realization of the DPA$(m,\beta)$ and let $\mathcal{G}^*(m)$ be a realization of its local weak limit. To start, for any $k \geq 1$ define
$$\mathbf{R}^{(k)}(n) = (1-c) \sum_{j=0}^k \left( \frac{c}{m} \right)^j \mathbf{1} A_n^j \qquad \text{and} \qquad  \mathcal{R}^{(k)} = (1-c) \left( 1 + \sum_{j=1}^k \left( \frac{c}{m} \right)^j \mathcal{P}_j^* \right) ,$$
where $A_n$ denotes the adjacency matrix of $\G_n$, and $\mathbf{R}^{(k)}(n)$, $\mathcal{R}^{(k)}$ are respectively computed on $\mathcal{G}_n$ and $\mathcal{G}^*(m)$. Observe that the scale-free PageRank on $\G_n$ is $\mathbf{R}(n) =\mathbf{R}^{(\infty)}(n)$, and $\mathcal{R} = \mathcal{R}^{(\infty)}$. Fix $\epsilon > 0$ and choose $\delta > 0$ and $k \geq 1$ such that
$$\pr{ \mathcal{R} \in ( r-\delta, r+\delta) } \leq \epsilon/6, \quad \delta^{-1} c^{k+1} \leq \epsilon/6 \quad \text{and} \quad  \pr{ \mathcal{R} - \mathcal{R}^{(k)} \geq \delta } \leq \epsilon/6,$$
which we can do since $r$ is a continuity point and $\mathcal{R}^{(k)} \to \mathcal{R}$ a.s.~as $k\to \infty$. 
Now let $V_n$ denote a uniformly chosen number in $\{1,\dots, n\}$ and note that 
$$\frac{1}{n} \sum_{i=1}^n \Identity(D_i^-(n) \geq l, R_i(n) > r)  = \pr{ \left. D_{V_n}^- \geq l, R_{V_n}(n) > r \right| \mathcal{G}_n }.$$
It follows that
\begin{align}
& \left| \pr{ \left. D_{V_n}^- \geq l, R_{V_n}(n) > r \right| \mathcal{G}_n } - \pr{ \mathcal{D}^- \geq l, \mathcal{R} > r} \right| \notag \\
&\leq \pr{ \left. D_{V_n}^- \geq l, R_{V_n}(n) > r \right| \mathcal{G}_n } -  \pr{ \left. D_{V_n}^- \geq l, R_{V_n}^{(k)}(n) > r \right| \mathcal{G}_n } \label{eq:G_nTail} \\
&\hspace{5mm} + \left|  \pr{ \left. D_{V_n}^- \geq l, R_{V_n}^{(k)}(n) > r \right| \mathcal{G}_n } - \pr{ \mathcal{D}^- \geq l, \mathcal{R}^{(k)} > r}  \right| \notag \\
&\hspace{5mm} + \pr{ \mathcal{D}^- \geq l, \mathcal{R} > r}  - \pr{ \mathcal{D}^- \geq l, \mathcal{R}^{(k)} > r} . \label{eq:LimitTail}
\end{align}
To bound \eqref{eq:G_nTail} note that it is equal to
\begin{align*}
&\pr{ \left. D_{V_n}^- \geq l, R_{V_n}(n) > r \geq R_{V_n}^{(k)}(n) \right| \mathcal{G}_n }  \\
&\leq \pr{ \left. R_{V_n}^{(k)}(n) + \delta > r \geq R_{V_n}^{(k)}(n) \right| \mathcal{G}_n } + \pr{ \left. R_{V_n}(n) - R_{V_n}^{(k)}(n) > \delta \right| \mathcal{G}_n } \\
&\leq \left|  \pr{ \left. R_{V_n}^{(k)}(n) \in (r-\delta, r] \right| \mathcal{G}_n } - \pr{ \mathcal{R}^{(k)} \in (r-\delta, r] } \right| \\
&\hspace{5mm} + \pr{ \mathcal{R} \in (r-\delta, r+\delta) } + \pr{ \mathcal{R} - \mathcal{R}^{(k)} \geq \delta} + \delta^{-1} \Expect{ \left. R_{V_n}(n) - R_{V_n}^{(k)}(n) \right| \mathcal{G}_n }.
\end{align*}
Now use the proof of Lemma~6.1 in \cite{olvera2019pagerank} to obtain that 
$$ \Expect{ \left. R_{V_n}(n) - R_{V_n}^{(k)}(n) \right| \mathcal{G}_n } = \frac{1}{n} \left\| \mathbf{R}(n) - \mathbf{R}^{(k)}(n) \right\|_1 \leq (1-c) \sum_{i=k+1}^\infty c^i =  c^{k+1}.$$
Hence, \eqref{eq:G_nTail} is bounded by
$$ \left|  \pr{ \left. R_{V_n}^{(k)}(n) \in (r-\delta, r] \right| \mathcal{G}_n } - \pr{ \mathcal{R}^{(k)} \in (r-\delta, r] } \right|  + 3\epsilon/6.$$
Similarly, \eqref{eq:LimitTail} is equal to
\begin{align*}
\pr{ \mathcal{D}^- \geq k, \mathcal{R} > r \geq \mathcal{R}^{(k)} } \leq \pr{ \mathcal{R} \in (r, r + \delta) } + \pr{ \mathcal{R} - \mathcal{R}^{(k)} \geq \delta } \leq 2\epsilon/6. 
\end{align*}
We have thus shown that
\begin{align}
&\pr{ \left| \pr{ \left. D_{V_n}^- \geq l, R_{V_n}(n) > r \right| \mathcal{G}_n } - \pr{ \mathcal{D}^- \geq l, \mathcal{R} > r} \right| > \epsilon } ] \notag \\
&\leq \pr{  \left|  \pr{ \left. R_{V_n}^{(k)}(n) \in (r-\delta, r] \right| \mathcal{G}_n } - \pr{ \mathcal{R}^{(k)} \in (r-\delta, r] } \right|  \right. \notag \\
&\hspace{5mm} \left. + \left|  \pr{ \left. D_{V_n}^- \geq l, R_{V_n}^{(k)}(n) > r \right| \mathcal{G}_n } - \pr{ \mathcal{D}^- \geq l, \mathcal{R}^{(k)} > r}  \right|  > \epsilon/6 }. \label{eq:UpperBd}
\end{align}
To see that this last probability converges to zero as $n \to \infty$, note that by Definition~3.6~(2) and Theorem~2.4~(3) in \cite{garavaglia2020local} we have that $\mathcal{G}_n$ converges to $\G^*(m)$ (with all vertices given the same mark $m$) in probability in the local weak sense. Further, observe that both the conditional probabilities in \eqref{eq:UpperBd} can be written in terms of expectations (with respect to the empirical measure $\mathbb{P}_n$ corresponding to selecting the root in $\G_n$ at random) of functions of the form:
$$
\Psi_{l,r,k}(G, \emptyset, M(G)) := \mathbb{I}\left(D_{\emptyset}^-(G, \emptyset, M(G)) \ge l, R^{(k)}_{\emptyset}(G, \emptyset, M(G)) >r\right).
$$
Here $(G,\emptyset,M(G))$ denotes a generic element of $\tilde{\mathbb{G}}_{\star}$ (the quotient space of directed, marked, rooted graphs equipped with a natural metric of directed local weak convergence), $D_v^-(G,\emptyset, M(G))$ denotes the in-degree of vertex $v$ in $(G,\emptyset, M(G))$, and $R^{(k)}_v(G, \emptyset, M(G))$ denotes its $k$-truncated PageRank. Since the function $\Psi_{l,r,k}$ is bounded and continuous with respect to the metric topology of $\tilde{\mathbb{G}}_{\star}$ for any fixed $l \in \mathbb{N}\cup \{0\}$, $k \in \mathbb{N}$ and $r > 0$, we conclude that \eqref{eq:UpperBd} converges to zero as $n \to \infty$. This completes the proof.
\end{proof}
Thus, quantifying the tail behavior of $\mathcal{D}^-$ and $\mathcal{R}$ provides information on the tail behavior of the in-degree and PageRank of a uniformly chosen vertex in $\G_n$ for large $n$. This is the main objective of this article.

\section{A unified description of $\G^*(m)$ for all $m \ge 1$}\label{loclimunif}

\
Section \ref{locsec} shows that the description of the local limit $\G^*(m)$ of the directed linear preferential attachment graph for $m \ge 2$ is starkly different from the $m=1$ case. However, we will see that the CTBP description of $\G^*(1)$ in terms of Yule processes is far more convenient for analyzing the distributional properties of PageRank, since one can use the explicit generators of Yule processes to build continuous time martingales. In particular, the convergence properties and moments of these martingales will turn out to be crucial tools for proving Theorems \ref{exponent} and \ref{rootpr}. Hence, it will be convenient to first obtain a new description of the P\'olya point graph $\G^*(m)$ for $m \ge 2$ in terms of Yule processes, which will enable a unified analysis of PageRank for all $m \geq 1$. 

\begin{theorem}\label{polyayule}
Fix $m \ge 1$. Consider the CTBP $\{\B^{(m-1 + \beta)}(t) : t \ge 0\}$ with root at $\emptyset$ where each individual $\oa \in \II$ reproduces independently according to $\xi_{\sigma_\oa}(t) = Y_\oa(t-\sigma_\oa), t \ge \sigma_\oa$, where $\{Y_\oa(\cdot) : \oa \in \II\}$ is an i.i.d. collection of $(m-1 + \beta)$-Yule Processes. Then
$$
\G^*(m) \equald \TT\left(\B^{(m-1 + \beta)}(\tau)\right)
$$
where $\tau$ is an $\operatorname{Exp}(2 +\beta/m)$ random variable that is independent of $\B^{(m-1+\beta)}(\cdot)$.
\end{theorem}

\begin{proof}
For $m=1$, the theorem follows directly from \cite[Theorem 1]{rudas-2}. For $m \ge 2$, the theorem follows from Proposition \ref{time} below.
\end{proof}

We will define and use the following notion of \emph{time-change} of a CTBP to prove the above theorem.

\begin{definition}\label{timechange}[Time-changed CTBP]
Consider two CTBPs $\{\B_{\sigma_\emptyset}(t) : t \ge \sigma_\emptyset \}$ and $\{\tilde{\B}_{\tilde{\sigma_\emptyset}}(t) : t \ge \tilde{\sigma_\emptyset}\}$ (possibly defined on different probability spaces) with roots $\emptyset, \tilde{\emptyset}$, born at times $\sigma_\emptyset, \tilde{\sigma_\emptyset}$ respectively, and associated reproduction point processes $\{\xi_x(t) : t \in [x, \infty)\}$, $x \in [\sigma_\emptyset, \infty),$ and $\{\tilde{\xi}_y(t) : t \in [y, \infty)\}$, $y \in [\tilde{\sigma_\emptyset}, \infty)$. For $\oa \in \II$, denote the birth times of $\oa$ in $\B_{\sigma_\emptyset}(\cdot)$ and $\tilde{\B}_{\tilde{\sigma_\emptyset}}(\cdot)$ by $\sigma_{\oa}$ and $\tilde{\sigma}_{\oa}$ respectively. For $T, \tilde{T}>0$, we say $\{\B_{\sigma_\emptyset}(s) : \sigma_\emptyset \le s \le T\}$ is a \emph{time-change} of $\{\tilde{\B}_{\tilde{\sigma_\emptyset}}(s) : \tilde{\sigma_\emptyset} \le s \le \tilde{T}\}$ if there exists a coupling $(\B_{\sigma_\emptyset}, \tilde{\B}_{\tilde{\sigma_\emptyset}})$ of the above CTBP such that for any $\oa \in \II$ with $\sigma_{\oa} \in [\sigma_\emptyset, T]$:
(i) $\tilde{\sigma}_{\oa} \in [\tilde{\sigma_\emptyset}, \tilde{T}]$ and
(ii) $\xi_{\sigma_{\oa}}(T) = \tilde{\xi}_{\tilde{\sigma}_{\oa}}(\tilde{T})$.

Note that (i) and (ii) make the time-change a symmetric relation. To see this, suppose $\{\B_{\sigma_\emptyset}(s) : \sigma_\emptyset \le s \le T\}$ is a time-change of $\{\tilde{\B}_{\tilde{\sigma_\emptyset}}(s) : \tilde{\sigma_\emptyset} \le s \le \tilde{T}\}$, and consider the above coupling satisfying (i) and (ii). Suppose there exists $\oa \in \II\setminus\{\emptyset\}$ such that $\tilde{\sigma}_{\oa} \in [\tilde{\sigma_\emptyset}, \tilde{T}]$ but $\sigma_{\oa} \notin [\sigma_\emptyset, T]$. Assume without loss of generality that $\oa$ is the minimal such element in $\II$. Write $\oa = (\mathbf{j},k)$ for some $\mathbf{j} \in \II$ and $k \in \mathbb{N}$. Then, the minimality of $\oa$ implies that $\sigma_{\mathbf{j}} \in [\sigma_\emptyset, T]$ and $\tilde{\sigma}_{\mathbf{j}} \in [\tilde{\sigma_\emptyset}, \tilde{T}]$. Thus, if  $\tilde{\sigma}_{\oa} \in [\tilde{\sigma_\emptyset}, \tilde{T}]$ but $\sigma_{\oa} \notin [\sigma_\emptyset, T]$, then $\xi_{\sigma_{\mathbf{j}}}(T) < \tilde{\xi}_{\tilde{\sigma}_{\mathbf{j}}}(\tilde{T})$, which is a contradiction to (ii). From this, we can deduce that $\{\tilde{\B}_{\tilde{\sigma_\emptyset}}(s) : \tilde{\sigma_\emptyset} \le s \le \tilde{T}\}$ is a time-change of  $\{\B_{\sigma_\emptyset}(s) : \sigma_\emptyset \le s \le T\}$.

In particular, under the time-change assertion, $\TT(\B_{\sigma_\emptyset}(T)) \equald \TT(\tilde{\B}_{\tilde{\sigma_\emptyset}}(\tilde{T}))$.
\end{definition}

By Remark \ref{polyactbp}, Theorem \ref{polyayule} immediately follows from the following proposition.

\begin{proposition}\label{time}
Fix $m \ge 1$. For $z>0$, let $\{\B^{pol}_{z}(t) : t \ge z\}$ be the CTBP from Remark \ref{polyactbp} and let $\{\B^{(m-1 + \beta)}(t) : t \ge 0\}$ be the CTBP from Theorem \ref{polyayule}. Then, there exists a coupling $(\sigma_\emptyset, \tau)$, where $\sigma_\emptyset \sim U^{\chi}$ for a $\operatorname{Uniform}(0,1)$ random variable $U$ and $\tau \sim \operatorname{Exp}(2 +\beta/m)$ such that, conditionally on $(\sigma_\emptyset, \tau)$, $\{\B^{pol}_{\sigma_\emptyset}(s) : \sigma_\emptyset \le s \le 1\}$ is a time-change of $\{\B^{(m-1+\beta)}(s) : 0 \le s \le \tau\}$.
\end{proposition}

The following lemma provides the key technical ingredient in proving Proposition~\ref{time}.

\begin{lemma}\label{pptoyule}
Let $\{T_i\}_{i \in \mathbb{N}}$ denote the arrival times of a Poisson point process on $\mathbb{R}_+$ with unit intensity. Let $\gamma \sim \operatorname{Gamma}(m+\beta, 1)$. Then the point process on $\mathbb{R}_+$ with arrival times given by
$$
S_i := \log\left(1+ \frac{T_i}{\gamma}\right), \ i \in \mathbb{N},
$$
is a $(m-1+\beta)$-Yule process.
\end{lemma}

To prove the lemma, we need the following fact that is easy to check.

\begin{lemma}\label{exptail}
Let $E \sim \operatorname{Exp}(1)$. For any $a,b,c>0$,
$$
\mathbb{E}\left(\mathbb{I}(E>ab)e^{-(a+E)c}\right) = \frac{e^{-a[c(1+b) + b]}}{c+1}.
$$
\end{lemma}

\begin{proof}[Proof of Lemma \ref{pptoyule}]
It suffices to show that for any $k \ge 2$,
\begin{equation}\label{prodexp}
(S_1, S_2-S_1, \dots, S_k - S_{k-1}) \sim \bigotimes_{i=1}^k \operatorname{Exp}(i+m-1+ \beta),
\end{equation}
where $\bigotimes$ denotes the product of independent random variables.
To show this, write $T_0=0, T_i = \sum_{s =1}^{i}E_{s}$, where $\{E_i : i \in \mathbb{N}\}$ are i.i.d. $\operatorname{Exp}(1)$ random variables. Observe that for any $t_1,\dots,t_k \in \mathbb{R}_+$,
\begin{multline}\label{pty1}
\pr{S_1 >t_1, S_2-S_1 > t_2, \dots, S_k - S_{k-1} > t_k}\\
= \pr{\frac{E_1}{\gamma} > e^{t_1} - 1, \frac{E_2}{\gamma + T_1} > e^{t_2} - 1, \dots, \frac{E_k}{\gamma + T_{k-1}} > e^{t_k} - 1}.
\end{multline}
We will now show that the following identity, denoted $I_j$, holds for any $1 \leq j \leq k-1$:
\begin{multline}\label{pty0}
 \pr{\frac{E_1}{\gamma} > e^{t_1} - 1, \frac{E_2}{\gamma + T_1} > e^{t_2} - 1, \dots, \frac{E_k}{\gamma + T_{k-1}} > e^{t_k} - 1}\\
 = \exp\left\lbrace -\sum_{\ell=0}^{j-1}(j-\ell -1) t_{k-\ell} \right\rbrace \mathbb{E}\left(\mathbb{I}\left(\frac{E_1}{\gamma} > e^{t_1} - 1, \dots, \frac{E_{k-j}}{\gamma + T_{k-j-1}} > e^{t_{k-j}} - 1\right)\right.\\
\left. \times \exp\left\lbrace -(\gamma + T_{k-j})(e^{\sum_{\ell=0}^{j-1}t_{k-\ell}} - 1)\right\rbrace\right).
\end{multline}
Our proof proceeds by induction in $j$.

%
For $j=1$, taking expectation over $E_k$ in the right hand side of \eqref{pty1} conditionally on $(\gamma, E_1,\dots, E_{k-1})$ yields
\begin{multline*}
\pr{\frac{E_1}{\gamma} > e^{t_1} - 1, \frac{E_2}{\gamma + T_1} > e^{t_2} - 1, \dots, \frac{E_k}{\gamma + T_{k-1}} > e^{t_k} - 1}\\
= \mathbb{E}\left(\mathbb{I}\left(\frac{E_1}{\gamma} > e^{t_1} - 1, \dots, \frac{E_{k-1}}{\gamma + T_{k-2}} > e^{t_{k-1}} - 1\right) e^{-[(\gamma + T_{k-1})(e^{t_k} - 1)]}\right)
\end{multline*}
establishing $I_1$. Suppose $k \ge 3$ and $I_j$ holds for some $1 \le j \le k-2$. Observe that
\begin{align}\label{pty2}
&\mathbb{E}\left(\mathbb{I}\left(\frac{E_1}{\gamma} > e^{t_1} - 1, \dots, \frac{E_{k-j}}{\gamma + T_{k-j-1}} > e^{t_{k-j}} - 1\right)\exp\left\lbrace -(\gamma + T_{k-j})(e^{\sum_{\ell=0}^{j-1}t_{k-\ell}} - 1)\right\rbrace\right)\\
&= \mathbb{E}\left(\mathbb{I}\left(\frac{E_1}{\gamma} > e^{t_1} - 1, \dots, \frac{E_{k-j-1}}{\gamma + T_{k-j-2}} > e^{t_{k-j-1}} - 1\right)\right.\notag\\
&\left. \times \mathbb{E}\left( \left. \mathbb{I}\left(\frac{E_{k-j}}{\gamma + T_{k-j-1}} > e^{t_{k-j}} - 1\right)\exp\left\lbrace -(\gamma + T_{k-j})(e^{\sum_{\ell=0}^{j-1}t_{k-\ell}} - 1)\right\rbrace \right| \gamma, E_1,\dots, E_{k-j-1}\right)\right).\notag
\end{align}
Applying Lemma \ref{exptail} with $E_{k-j}$ in place of $E$, $a= \gamma + T_{k-j-1}$, $b = e^{t_{k-j}}-1$ and $c= e^{\sum_{\ell=0}^{j-1}t_{k-\ell}} - 1$, we obtain
\begin{multline}\label{pty3}
\mathbb{E}\left( \left. \mathbb{I}\left(\frac{E_{k-j}}{\gamma + T_{k-j-1}} > e^{t_{k-j}} - 1\right)\exp\left\lbrace -(\gamma + T_{k-j})(e^{\sum_{\ell=0}^{j-1}t_{k-\ell}} - 1)\right\rbrace \right| \gamma, E_1,\dots, E_{k-j-1}\right)\\
= \exp\left\lbrace -\sum_{\ell=0}^{j-1}t_{k-\ell}\right\rbrace \exp\left\lbrace -(\gamma + T_{k-j-1})\left(e^{\sum_{\ell=0}^{j}t_{k-\ell}}-1\right)\right\rbrace.
\end{multline}
Substituting \eqref{pty3} into \eqref{pty2} gives
\begin{multline}\label{pty4}
\mathbb{E}\left(\mathbb{I}\left(\frac{E_1}{\gamma} > e^{t_1} - 1, \dots, \frac{E_{k-j}}{\gamma + T_{k-j-1}} > e^{t_{k-j}} - 1\right)\exp\left\lbrace -(\gamma + T_{k-j})(e^{\sum_{\ell=0}^{j-1}t_{k-\ell}} - 1)\right\rbrace\right)\\
= \exp\left\lbrace -\sum_{\ell=0}^{j-1}t_{k-\ell}\right\rbrace \mathbb{E}\left(\mathbb{I}\left(\frac{E_1}{\gamma} > e^{t_1} - 1, \dots, \frac{E_{k-j-1}}{\gamma + T_{k-j-2}} > e^{t_{k-j-1}} - 1\right) \right.\\\left. \times \exp\left\lbrace -(\gamma + T_{k-j-1})\left(e^{\sum_{\ell=0}^{j}t_{k-\ell}}-1\right)\right\rbrace \right).
\end{multline}
As $I_j$ holds, using \eqref{pty4} in \eqref{pty0} gives,
\begin{align*}
&\pr{\frac{E_1}{\gamma} > e^{t_1} - 1, \frac{E_2}{\gamma + T_1} > e^{t_2} - 1, \dots, \frac{E_k}{\gamma + T_{k-1}} > e^{t_k} - 1}\\
&= \exp\left\lbrace -\sum_{\ell=0}^{j-1}(j-\ell -1) t_{k-\ell} \right\rbrace \exp\left\lbrace -\sum_{\ell=0}^{j-1}t_{k-\ell}\right\rbrace\\
& \hspace{30pt}\times \mathbb{E}\left(\mathbb{I}\left(\frac{E_1}{\gamma} > e^{t_1} - 1, \dots, \frac{E_{k-j-1}}{\gamma + T_{k-j-2}} > e^{t_{k-j-1}} - 1\right)\right.\\
&\hspace{100pt}\left. \times\exp\left\lbrace -(\gamma + T_{k-j-1})\left(e^{\sum_{\ell=0}^{j}t_{k-\ell}}-1\right)\right\rbrace \right)\\
&=\exp\left\lbrace -\sum_{\ell=0}^{j}(j-\ell) t_{k-\ell} \right\rbrace\\
&\hspace{30pt} \times \mathbb{E}\left(\mathbb{I}\left(\frac{E_1}{\gamma} > e^{t_1} - 1, \dots, \frac{E_{k-j-1}}{\gamma + T_{k-j-2}} > e^{t_{k-j-1}} - 1\right)\right.\\
&\hspace{100pt}\left. \times \exp\left\lbrace -(\gamma + T_{k-j-1})\left(e^{\sum_{\ell=0}^{j}t_{k-\ell}}-1\right)\right\rbrace \right).
\end{align*}
This proves $I_{j+1}$. Thus, by induction, for any $k \ge 2$, \eqref{pty0} holds for all $1 \le j \le k-1$. 

To complete the proof of the lemma, note that from \eqref{pty0} with $j=k-1$, we obtain
\begin{align}\label{pty5}
&\pr{\frac{E_1}{\gamma} > e^{t_1} - 1, \frac{E_2}{\gamma + T_1} > e^{t_2} - 1, \dots, \frac{E_k}{\gamma + T_{k-1}} > e^{t_k} - 1}\\
 &= \exp\left\lbrace -\sum_{\ell=0}^{k-2}(k-\ell -2) t_{k-\ell} \right\rbrace \mathbb{E}\left(\mathbb{I}\left(\frac{E_1}{\gamma} > e^{t_1} - 1\right)\exp\left\lbrace -(\gamma + E_{1})(e^{\sum_{\ell=0}^{k-2}t_{k-\ell}} - 1)\right\rbrace\right)\notag\\ 
 &=  \exp\left\lbrace -\sum_{\ell=0}^{k-2}(k-\ell -2) t_{k-\ell}\right\rbrace \exp\left\lbrace -\sum_{\ell=0}^{k-2}t_{k-\ell}\right\rbrace \mathbb{E}\left(\exp\left\lbrace -\gamma\left(e^{\sum_{\ell=0}^{k-1}t_{k-\ell}}-1\right)\right\rbrace\right)\notag\\
 &=  \exp\left\lbrace -\sum_{\ell=0}^{k-1}(k-\ell -1) t_{k-\ell}\right\rbrace \mathbb{E}\left(\exp\left\lbrace -\gamma\left(e^{\sum_{\ell=0}^{k-1}t_{k-\ell}}-1\right)\right\rbrace\right), \notag
\end{align}
where the second equality follows from Lemma \ref{exptail} with $E_{1}$ in place of $E$, $a= \gamma$, $b = e^{t_{1}}-1$ and $c= e^{\sum_{\ell=0}^{k-2}t_{k-\ell}} - 1$. As $\gamma \sim \operatorname{Gamma}(m+\beta, 1)$, the explicit form of the moment generating function for Gamma distributions gives
$$
\mathbb{E}\left(\exp\left\lbrace -\gamma\left(e^{\sum_{\ell=0}^{k-1}t_{k-\ell}}-1\right)\right\rbrace\right) = e^{-(m+\beta)\sum_{\ell=0}^{k-1}t_{k-\ell}}.
$$
Using this in \eqref{pty5}, we have
\begin{equation*}
\pr{\frac{E_1}{\gamma} > e^{t_1} - 1, \frac{E_2}{\gamma + T_1} > e^{t_2} - 1, \dots, \frac{E_k}{\gamma + T_{k-1}} > e^{t_k} - 1} = \exp\left\lbrace -\sum_{i=1}^{k}(i+m-1+\beta) t_i\right\rbrace.
\end{equation*}
Thus, from \eqref{pty1}, we conclude \eqref{prodexp} which proves the lemma.
\end{proof}

We are now ready to prove Proposition~\ref{time}.

\begin{proof}[Proof of Proposition \ref{time}]
Define the function $\tau: (0,1) \rightarrow \mathbb{R}_+$ by $\tau(z) := \frac{m}{m+\beta}\log\left(\frac{1}{z}\right)$. Sample $\sigma_\emptyset \sim U^{\chi}$ for a $\operatorname{Uniform}(0,1)$ random variable $U$ and set $\tau := \tau(\sigma_\emptyset)$. Then, it follows that $\tau \sim \operatorname{Exp}(2 +\beta/m)$. This gives a coupling $(\sigma_\emptyset, \tau)$ required by the proposition. To prove the proposition, it suffices to show that, for any $z>0$, $\{\B^{pol}_{z}(s) : z \le s \le 1\}$ is a time-change of $\{\B^{(m-1+\beta)}(s) : 0 \le s \le \tau(z)\}$.

Fix $z>0$. Recall the CTBP $\{\B^{pol}_{z}(t) : t \ge z\}$ with associated reproduction point process for each $x \in [z, \infty)$ given by the  \textcolor{black}{nonhomogeneous} Poisson point process $\{\xi^{pol}_x(t) : t \in [x, \infty)\}$ with (random) intensity $\rho_x(v) := \gamma_{x}\frac{\psi v^{\psi-1}}{x^{\psi}}, \ v \in [x,\infty)$. For $\oa \in \II$, denote the birth time of individual $\oa$ into $\B^{pol}_{z}(\cdot)$ by $\sigma_{\oa}$. For $0 \le y < \infty$, define the point process
\begin{equation}\label{coupdef}
\tilde{\xi}_y(s) := \xi^{pol}_{ze^{(1+\beta/m)y}}(ze^{(1+\beta/m)s}), \ s \ge y.
\end{equation}
Consider the CTBP $\{\tilde{\B}(t) : t \ge 0\}$ obtained from the reproduction point processes defined in \eqref{coupdef}. We suppress the dependence of the CTBP on $z$ for notational convenience. Denote the birth time of individual $\oa \in \II$ into this CTBP by $\tilde{\sigma}_{\oa}$.

To show the time-change assertion, it suffices to show the following:
\begin{itemize}
\item[(I)] For each $\oa \in \II$, $Y_{\oa}(t) := \tilde{\xi}_{\tilde{\sigma}_{\oa}}(t + \tilde{\sigma}_{\oa}), t \ge 0,$ is a $(m-1+\beta)$-Yule process; and
\item[(II)] for any $\oa \in \II$ such that $\sigma_{\oa} \in [z,1]$, $\tilde{\sigma}_{\oa} \in [0, \tau(z)]$ and $\tilde{\xi}_{\tilde{\sigma}_{\oa}}(\tau(z)) = \xi^{pol}_{\sigma_{\oa}}(1)$.
\end{itemize}
To prove (I), note that for any $\oa \in \II$, writing $w := ze^{(1 + \beta/m)\tilde{\sigma}_{\oa}}$ gives,
\begin{equation}\label{p1}
Y_{\oa}(t) := \tilde{\xi}_{\tilde{\sigma}_{\oa}}(t + \tilde{\sigma}_{\oa}) =  \xi^{pol}_{ze^{(1+\beta/m)\tilde{\sigma}_{\oa}}}(ze^{(1+\beta/m)(t + \tilde{\sigma}_{\oa})}) = \xi^{pol}_{w}(w e^{(1+\beta/m)t}).
\end{equation}
As $\xi^{pol}_{w}(\cdot)$ is a \textcolor{black}{nonhomogeneous} Poisson point process with intensity given by $\rho_w(v) := \gamma_w\frac{\psi v^{\psi - 1}}{w^{\psi}}, v \in [w,\infty)$, we can write it in terms of a Poisson point process $N(\cdot)$ on $\mathbb{R}_+$ of unit intensity (independent of $\gamma_w$) in the form
\begin{equation*}
\xi^{pol}_{w}(t) = N\left(\gamma_w\left[\left(t/w\right)^{\psi} - 1\right]\right), \ t \ge w.
\end{equation*}
Thus, the arrival times $\{S'_i\}_{i \in \mathbb{N}}$ of $\xi^{pol}_{w}(\cdot)$ can be expressed in terms of arrival times $\{T_i\}_{i \in \mathbb{N}}$ of $N(\cdot)$ as
$$
S'_i = w\left(1 + \frac{T_i}{\gamma_w}\right)^{1/\psi}, \ i \in \mathbb{N}.
$$
By \eqref{p1}, the arrival times $\{S_i\}_{i \in \mathbb{N}}$ of $Y_{\oa}(\cdot)$ can be expressed as
$$
S_i = \frac{m}{m+\beta}\log\left(\frac{S'_i}{w}\right) = \log\left(1 + \frac{T_i}{\gamma_w}\right), \ i \in \mathbb{N}.
$$
Thus, by Lemma \ref{pptoyule}, $Y_{\oa}(\cdot)$ is a $(m-1+\beta)$-Yule process, proving (I).

To show (II), we will establish the following:
\begin{equation}\label{birth}
\tilde{\sigma}_{\oa} = \frac{m}{m + \beta}\log \left(\frac{\sigma_{\oa}}{z}\right) \ \text{ for every } \oa \in \II.
\end{equation}
We will prove \eqref{birth} by induction on the length of $\oa$. Clearly, \eqref{birth} is true for $\oa = \emptyset$. Suppose \eqref{birth} is true for all $\oa \in \mathbb{N}^j$ for some $j \ge 0$. Take any $\oa \in \mathbb{N}^{j+1}$. Then there is $\oa' \in \mathbb{N}^j$ and $\ell \in \mathbb{N}$ such that $\oa = (\oa', \ell)$. Note that, by induction hypothesis, 
$
ze^{(1+\beta/m)\tilde{\sigma}_{\oa'}} = \sigma_{\oa'}.
$
Thus,
$$
\tilde{\xi}_{\tilde{\sigma}_{\oa'}}(s) := \xi^{pol}_{ze^{(1+\beta/m)\tilde{\sigma}_{\oa'}}}(ze^{(1+\beta/m)s}) = \xi^{pol}_{\sigma_{\oa'}}(ze^{(1+\beta/m)s}), \ s \ge \tilde{\sigma}_{\oa'}.
$$
Hence,
$$
\sigma_{\oa} = \sigma_{(\oa',\ell)} = ze^{(1+\beta/m)\tilde{\sigma}_{(\oa',\ell)}} = ze^{(1+\beta/m)\tilde{\sigma}_{\oa}}
$$
proving \eqref{birth} for all $\oa \in \mathbb{N}^{j+1}$. By induction \eqref{birth} holds for all $\oa \in \II$.
By \eqref{birth}, for any $\oa \in \II$,
\begin{equation}\label{p2}
\sigma_{\oa} \in [z,1] \ \text{ if and only if } \tilde{\sigma}_{\oa} \in [0, \tau(z)],
\end{equation}
and, recalling $\tau(z) = \frac{m}{m+\beta}\log\left(\frac{1}{z}\right)$,
\begin{equation}\label{p3}
\tilde{\xi}_{\tilde{\sigma}_{\oa}}(\tau(z)) = \xi^{pol}_{ze^{(1+\beta/m)\tilde{\sigma}_{\oa}}}(ze^{(1+\beta/m)\tau(z)}) = \xi^{pol}_{\sigma_{\oa}}(ze^{(1+\beta/m)\tau(z)}) = \xi^{pol}_{\sigma_{\oa}}(1).
\end{equation}
(II) follows from \eqref{p2} and \eqref{p3}. This completes the proof of the proposition.
\end{proof}

\section{Proofs: Power law exponent of PageRank} \label{proofs}
In this section, we will prove Theorem~\ref{exponent}.

We will work with fixed $m, \beta$ and abbreviate $\B^{(m-1+\beta)}(\cdot)$ as $\B(\cdot)$. We will also write $\G(t) := \TT(\B(t)), t \ge 0$. For $i \ge 1, t \ge 0$, let $\PP_i(t)$ denote the number of (directed) paths in $\G(t)$ of length $i$ that end at the root. By \eqref{limpr} and Theorem \ref{polyayule}, note that
\begin{equation}\label{prdist}
\mathcal{R} \equald (1-c)\left(1 + \sum_{k=1}^{\infty}\left(\frac{c}{m}\right)^k \PP_k(\tau)\right), \ \ \tau \sim \operatorname{Exp}(2+ \beta/m).
\end{equation}
Thus, distributional analysis of $\mathcal{R}$ reduces to that of the right hand side of equation \eqref{prdist}.

Observe that
$
\boldsymbol{\mathcal{P}}(t) := (\PP_1(t), \PP_2(t), \dots)^T, \ t \ge 0,
$
is a $\mathbb{N}^{\infty}$-valued Markov process whose generator is given by
\begin{align}\label{gendef}
\LL f(\mathbf{p}) &:= (f(\mathbf{p} + \mathbf{e}^{(1)}) - f(\mathbf{p})) (p_1 + m + \beta)\\
&\quad + \sum_{i=2}^{\infty}(f(\mathbf{p} + \mathbf{e}^{(i)}) - f(\mathbf{p})) (p_i + (m + \beta)p_{i-1}), \ \ \mathbf{p} = (p_1,p_2,\dots) \in \mathbb{N}^{\infty},\notag
\end{align}
for any function $f : \mathbb{N}^{\infty} \rightarrow \mathbb{R}$ for which the right hand side above is well-defined, and $\mathbf{e}^{(i)}$ denotes the vector with a one in the $it$h coordinate and zeros elsewhere. Thus, writing $\LL \boldsymbol{\mathcal{P}}(\cdot) := (\LL \PP_1(\cdot), \LL \PP_2(\cdot),\dots)^T$, we get the following representation:
\begin{equation}\label{jordan}
\LL \boldsymbol{\mathcal{P}}(t) = \textcolor{black}{Q}\boldsymbol{\mathcal{P}}(t) + (m+\beta) \mathbf{e}^{(1)}, \ t \ge 0,black
\end{equation}
where \textcolor{black}{$Q$} is the infinite Jordan-block-type matrix given by \textcolor{black}{$Q_{ii} = 1$} for $i \ge 1$,  \textcolor{black}{$Q_{i (i-1)} = m+ \beta$} for $i \ge 2$ and  \textcolor{black}{$Q_{ij} = 0$} otherwise. Note that the representation is always well-defined as for any $t \ge 0$, there exists $i=i(t) \in \mathbb{N}$ such that $\PP_j(t) = 0$ for all $j > i$.

For $n \in \mathbb{N}$, write the $n \times n$ principal sub-matrix of \textcolor{black}{$Q$} as \textcolor{black}{$Q_n$}, and let $\boldsymbol{\mathcal{P}}_n(\cdot) := (\PP_1(\cdot), \dots, \PP_n(\cdot))^T$. The Jordan-block-type form of \textcolor{black}{$Q_n$} implies for any $t \in \mathbb{R}$ \cite[Chapter 2, Theorem 1.4]{weintraub2008jordan},
$$
(e^{\textcolor{black}{Q_n} t})_{ij} =
\left\{
	\begin{array}{ll}
		\frac{(m+\beta)^{i-j}t^{i-j}}{(i-j)!}e^t   & \mbox{if } i \ge j, \\
		0 & \mbox{if } i < j. 
	\end{array}
\right.
$$
Define the $\mathbb{R}^n$-valued process
$$
\mathbf{M}_n(t) := e^{-\textcolor{black}{Q_n} t}\boldsymbol{\mathcal{P}}_n(t) - (m+\beta)\int_0^t e^{-\textcolor{black}{Q_n} s} \mathbf{e}_n^{(1)}ds, \quad t \ge 0.
$$
For $1 \le i \le n$, we will write $M_{i}(\cdot)$ for the $i$-th coordinate of $\mathbf{M}_n(\cdot)$. Note that due to the lower triangular form of \textcolor{black}{$Q_n$}, $M_i(\cdot)$ does not depend on $n$.

\begin{lemma}\label{mn}
For any $n \in \mathbb{N}$, $\mathbf{M}_n(\cdot)$ is a martingale with respect to the natural filtration generated by $\boldsymbol{\mathcal{P}}(\cdot)$.
\end{lemma}

\begin{proof}
Fix $n \in \mathbb{N}$. For any $t \ge 0$,
\begin{align*}
\mathbb{E}\left(\|\mathbf{M}_n(t)\|_1\right) &\le \mathbb{E}\left( \sum_{i=1}^n  \left|\sum_{j=1}^n (e^{-\textcolor{black}{Q_n} t})_{ij}  \mathcal{P}_j(t)\right|  \right) + (m+\beta)\int_0^t \mathbf{1} e^{-\textcolor{black}{Q_n} s} \mathbf{e}_n^{(1)}ds \\
&\le (m+\beta)^{n-1}\mathbb{E}\left(|\B(t)|\right) + (m+\beta)\int_0^t \mathbf{1} e^{-\textcolor{black}{Q_n} s} \mathbf{e}_n^{(1)}ds
\end{align*}
where $|\B(t)| = 1+ \sum_{i=1}^{\infty}\PP_i(t)$ denotes the total number of individuals in the CTBP $\B(\cdot)$ at time $t$, and $\mathbf{1} := (1,1,\dots,1)$. As $m+1+\beta$ is the Malthusian rate for the CTBP $\B(\cdot)$ (see \cite[Section~4.2]{rudas-2}), by \cite[Proposition~2.2]{nerman1981convergence}, 
$$\sup_{t \ge 0}e^{-(m+1+\beta)t}\mathbb{E}\left(|\B(t)|\right) < \infty.$$ 
Hence, from the above bound, $\mathbb{E}\left(\|\mathbf{M}_n(t)\|_1\right) < \infty$ for all $t \ge 0$. 

Thus, to show that $\mathbf{M}_n(\cdot)$ is a martingale, it remains to show that black
$$
\LL\mathbf{M}_n(t) := (\LL M_{1}(t),\dots, \LL M_{n}(t))^T = \mathbf{0} \ \text{ for all } \ t \ge 0,
$$
where $\LL M_i(t) := \frac{d}{ds} \mathbb{E}(M_i(s)) \mid_{s=t}$, $1 \le i \le n$. This follows from the calculation
\begin{align*}
\LL\mathbf{M}_n(t) &= -\textcolor{black}{Q_n} e^{-\textcolor{black}{Q_n} t} \boldsymbol{\mathcal{P}}_n(t) + e^{-\textcolor{black}{Q_n} t} \LL \boldsymbol{\mathcal{P}}_n(t) - (m+\beta) e^{-\textcolor{black}{Q_n} t} \mathbf{e}_n^{(1)}\\
&= -\textcolor{black}{Q_n} e^{-\textcolor{black}{Q_n} t} \boldsymbol{\mathcal{P}}_n(t) + e^{-\textcolor{black}{Q_n} t}\left(\textcolor{black}{Q_n} \boldsymbol{\mathcal{P}}_n(t) + (m+\beta) \mathbf{e}_n^{(1)}\right) -(m+\beta) e^{-\textcolor{black}{Q_n} t} \mathbf{e}_n^{(1)}=0.
\end{align*}
This proves the lemma.
\end{proof}

For $t \in \mathbb{R}$, define the exponential $e^{\textcolor{black}{Q} t}$ of the infinite matrix \textcolor{black}{$Q$} by requiring the $n\times n$ principal sub-matrix of $e^{\textcolor{black}{Q} t}$ to equal $e^{\textcolor{black}{Q_n}t}$ for any $n \times n$. Again, the lower triangular form of $\textcolor{black}{Q_n}$, and hence $e^{\textcolor{black}{Q_n}t}$, implies that $e^{\textcolor{black}{Q} t}$ is well-defined. Define the $\mathbb{R}^{\infty}$-valued process
\begin{equation}\label{mdef}
\mathbf{M}(t) := e^{-\textcolor{black}{Q} t}\boldsymbol{\mathcal{P}}(t) - (m+\beta)\int_0^t e^{-\textcolor{black}{Q} s} \mathbf{e}^{(1)}ds =: e^{-\textcolor{black}{Q} t} \boldsymbol{\mathcal{P}} (t) - \boldsymbol{\phi}(t), \ t \ge 0.
\end{equation}
By Lemma \ref{mn}, the coordinate processes of $\mathbf{M}(\cdot)$ are martingales. Define
\begin{equation}\label{ptdef}
\PP^*(t) := \sum_{k=1}^\infty \left(\frac{c}{m}\right)^k\PP_k(t), \ t \ge 0.
\end{equation} 
This can be thought of as the `continuous time PageRank' as $\PP^*(\tau) = (1-c)^{-1}\mathcal{R}-1$ for an independently sampled $\tau \sim \operatorname{Exp}(2+ \beta/m)$. By defining the vector $\mathbf{v} := ((c/m),(c/m)^2,\dots)$ and using \eqref{mdef}, we obtain that
\begin{equation}\label{prm}
\PP^*(t) = \mathbf{v}\left[e^{\textcolor{black}{Q} t}\left(\mathbf{M}(t) + \boldsymbol{\phi}(t)\right)\right], \ t \ge 0.
\end{equation}
Note that for any $i \in \mathbb{N}$, $t \in \mathbb{R}$,
\begin{align}\label{cint}
\left(\mathbf{v} e^{\textcolor{black}{Q} t}\right)_i &= \sum_{j=i}^{\infty}\left(\frac{c}{m}\right)^j\left(e^{\textcolor{black}{Q} t}\right)_{ji} = \sum_{j=i}^{\infty}\left(\frac{c}{m}\right)^j\frac{(m+\beta)^{j-i}t^{j-i}}{(j-i)!}e^t\\
&= \left(\frac{c}{m}\right)^{i}e^t  \sum_{j=0}^{\infty}\left(\frac{c}{m}\right)^j\frac{(m+\beta)^{j}t^{j}}{j!} = \left(\frac{c}{m}\right)^i e^{(1+ (m+\beta)c/m)t}.\notag
\end{align}
Hence, by \eqref{prm}, writing $\boldsymbol{\phi}(\cdot) = (\phi_1(\cdot), \phi_2(\cdot),\dots)^T$, and
$
\theta := 1+ (m+\beta)c/m,
$
\begin{align}\label{martdec}
\PP^*(t) &= \left[\mathbf{v} e^{\textcolor{black}{Q} t}\right]\left(\mathbf{M}(t) + \boldsymbol{\phi}(t)\right) =  e^{(1+ (m+\beta)c/m)t}\sum_{i=1}^{\infty} \left(\frac{c}{m}\right)^i \left(M_i(t) + \phi_i(t)\right) \\
&=  e^{\theta t}\left(\sum_{i=1}^{\infty} \left(\frac{c}{m}\right)^i M_i(t) + \sum_{i=1}^{\infty} \left(\frac{c}{m}\right)^i\phi_i(t)\right) =: e^{\theta t}(M^*(t) + \phi^*(t)),   \quad t\ge 0.\notag
\end{align}
To justify the first and third equalities above, we will verify that the sums $\sum_{i=1}^{\infty} \left(\frac{c}{m}\right)^i |M_i(t)|$ and $\sum_{i=1}^{\infty}\left(\frac{c}{m}\right)^i|\phi_i(t)|$ are almost surely finite. This is done as follows. Note that, writing $\overline{\mathbf{M}}(\cdot) := (|M_1(\cdot)|, |M_2(\cdot)|, \dots)^T$,
$$
\overline{\mathbf{M}}(t) \le e^{\textcolor{black}{Q} t}\boldsymbol{\mathcal{P}}(t) + (m+\beta)\int_0^t e^{\textcolor{black}{Q} s} \mathbf{e}^{(1)}ds.
$$
Thus, using \eqref{cint},
\begin{align}\label{martfin}
\sum_{i=1}^{\infty} \left(\frac{c}{m}\right)^i |M_i(t)| &\le e^{\theta t}\sum_{i=1}^{\infty} \left(\frac{c}{m}\right)^i \PP_i(t) + (m + \beta)\left(\frac{c}{m}\right)\int_0^t e^{\theta s} ds\notag\\
&\le e^{\theta t} \left(\frac{c}{m}\right) |\B(t)| + (m + \beta)\left(\frac{c}{m}\right)\int_0^t e^{\theta s} ds.
\end{align}
As $\sup_{t \ge 0}e^{-(m+1+\beta)t}\mathbb{E}\left(|\B(t)|\right) < \infty$, the above sum has finite expected value, and hence, is almost surely finite. Similarly,
\begin{equation}\label{phifin}
\sum_{i=1}^{\infty}\left(\frac{c}{m}\right)^i|\phi_i(t)| \le (m + \beta)\left(\frac{c}{m}\right)\int_0^t e^{\theta s} ds.
\end{equation}
The above bounds validate \eqref{martdec} and also show that $M^*(\cdot)$ and $\phi^*(\cdot)$ are almost surely well-defined and finite-valued.

The following lemma is a crucial technical ingredient in proving Theorem \ref{exponent}.

\begin{lemma}\label{lpbdd}
$M^*(\cdot)$ is a martingale with respect to the natural filtration generated by $\boldsymbol{\mathcal{P}}(\cdot)$. Moreover, for any $k \in \mathbb{N}$,
\begin{equation}\label{lp}
\sup_{t \ge 0} \mathbb{E}\left(|M^*(t)|^k\right) \le 2^k\left[\left(\frac{(m+\beta)c}{\theta m}\right)^k + \left(\frac{m+\beta + \theta}{\theta}\right)^k 2^{(k+1)(k+2)/2}\right].
\end{equation}
\end{lemma}

\begin{proof}
By Lemma \ref{mn}, $M_n^*(\cdot) := \sum_{i=1}^{n} \left(\frac{c}{m}\right)^i M_i(\cdot)$ is a martingale for each $n \in \mathbb{N}$. Moreover, the steps leading to \eqref{martfin} also give for each $t \ge 0$,
$
\lim_{n \rightarrow \infty} \mathbb{E}\left(|M^*(t) - M^*_n(t)|\right) = 0.
$
This implies that $M^*(\cdot)$ is a martingale.

To prove \eqref{lp}, note that
\begin{equation}\label{phibd}
\phi^*(t) = (m+\beta)\int_0^t[\mathbf{v} e^{-\textcolor{black}{Q} s}]\mathbf{e}^{(1)}ds = (m+\beta)\left(\frac{c}{m}\right)\int_0^{t}e^{-\theta s} ds \le \frac{(m+\beta)c}{\theta m},
\end{equation}
where the first equality uses \eqref{mdef} and \eqref{martdec}, and the second equality follows from \eqref{cint}. In particular, the second equality shows that $\phi^*(\cdot)$ is non-negative and non-decreasing. Hence, from \eqref{martdec},
$$
|M^*(t)| \le e^{-\theta t} \PP^*(t) + \frac{(m+\beta)c}{\theta m}.
$$
Write $U(t) :=  e^{-\theta t} \PP^*(t)$, $t \ge 0$. By the above bound, to prove \eqref{lp}, it suffices to show that
\begin{equation}\label{ubd}
\sup_{t \ge 0} \mathbb{E}\left(U(t)^\ell\right) \le \left(\frac{m+\beta + \theta}{\theta}\right)^\ell 2^{(\ell+1)(\ell+2)/2} =: D_\ell, \quad \ell \ge 1.
\end{equation}
We will show \eqref{ubd} by induction on $\ell$. For $\ell=1$ and any $t \ge 0$, note that using \eqref{phibd} and the fact that $M^*(\cdot)$ is a martingale, and therefore  $\mathbb{E}(M^*(t)) = M^*(0) = 0$, we have
$$
\mathbb{E}(U(t)) = \mathbb{E}(M^*(t)) + \phi^*(t) \le \frac{(m+\beta)c}{\theta m} \le D_1.
$$
Suppose for some $k \ge 2$, \eqref{ubd} holds for all $\ell \le k-1$. Let $f_k: \mathbb{N}^{\infty} \rightarrow \mathbb{R}$ be defined by $f_k(\mathbf{p}) = \left(\sum_{i=1}^{\infty} (c/m)^i p_i\right)^k$. Applying the generator $\LL$ to $f_k$ by the formula \eqref{gendef} and setting $\PP_0(\cdot) =0$, we obtain
\begin{equation}\label{udiff}
\left. \frac{d}{ds}\mathbb{E}\left(U(s)^k\right) \right|_{s=t} \ = \mathbb{E}\left(\Lambda(t)\right) ,
\end{equation}
where
\begin{align}\label{b1}
\Lambda(t) &= - k\theta U(t)^k + e^{-\theta k t}\LL f_k(\boldsymbol{\mathcal{P}}(t))\\
&= - k\theta U(t)^k + e^{-\theta k t} \left((\PP^*(t) + c/m)^k - \PP^*(t)^k\right)(\PP_1(t) + m+ \beta)\notag\\
&\quad + e^{-\theta k t} \sum_{j=2}^{\infty}\left((\PP^*(t) + (c/m)^j)^k - \PP^*(t)^k\right)(\PP_j(t) + (m+ \beta)\PP_{j-1}(t))\notag\\
&= - k\theta U(t)^k + e^{-\theta k t} \left((\PP^*(t) + c/m)^k - \PP^*(t)^k\right)(m+ \beta)\notag\\
&\quad + e^{-\theta k t} \sum_{j=1}^{\infty}\left(\sum_{i=0}^{k-1}{k\choose i}\PP^*(t)^i(c/m)^{j(k-i)}\right)(\PP_j(t) + (m+ \beta)\PP_{j-1}(t)).\notag
\end{align}
As $c \in (0,1)$,
\begin{align*}
&\sum_{j=1}^{\infty}\left(\sum_{i=0}^{k-1}{k\choose i}\PP^*(t)^i(c/m)^{j(k-i)}\right)(\PP_j(t) + (m+ \beta)\PP_{j-1}(t))\\
&\le \sum_{i=0}^{k-1}{k\choose i}\PP^*(t)^i\sum_{j=1}^{\infty}(c/m)^j\PP_j(t) + \frac{(m+\beta) c}{m}\sum_{i=0}^{k-1}{k\choose i}\PP^*(t)^i\sum_{j=1}^{\infty}(c/m)^{j-1}\PP_{j-1}(t)\\
&= \left( 1 + \frac{(m+\beta)c}{m} \right) \sum_{i=0}^{k-1}{k\choose i}\PP^*(t)^{i+1} = \theta \sum_{i=0}^{k-1}{k\choose i}\PP^*(t)^{i+1}.
\end{align*}
Using this bound in \eqref{b1}, we obtain
\begin{align*}
\Lambda(t) &\le - k\theta U(t)^k + e^{-\theta k t} \left((\PP^*(t) + c/m)^k - \PP^*(t)^k\right)(m+ \beta) + \theta e^{-\theta k t} \sum_{i=0}^{k-1}{k\choose i}\PP^*(t)^{i+1}\\
&= e^{-\theta k t} \left((\PP^*(t) + c/m)^k - \PP^*(t)^k\right)(m+ \beta) + \theta e^{-\theta k t} \sum_{i=0}^{k-2}{k\choose i}\PP^*(t)^{i+1}\\
&= (m+\beta) e^{-\theta k t}\sum_{i=0}^{k-1}{k\choose i}\PP^*(t)^{i}(c/m)^{k-i} + \theta e^{-\theta k t} \sum_{i=0}^{k-2}{k\choose i}\PP^*(t)^{i+1}\\
&\le (m+\beta) e^{-\theta k t}\sum_{i=0}^{k-1}{k\choose i}\PP^*(t)^{i} + \theta e^{-\theta k t} \sum_{i=1}^{k-1}{k\choose i-1}\PP^*(t)^{i}\\
&\le (m+\beta) e^{-\theta t}\sum_{i=0}^{k-1}{k\choose i}U(t)^{i} + \theta e^{-\theta t} \sum_{i=1}^{k-1}{k\choose i-1}U(t)^{i}\\
&\le (m+\beta + \theta)e^{-\theta t} \sum_{i=0}^{k-1}{k+ 1 \choose i}U(t)^{i}.
\end{align*}
For the first equality in the above display, we have used $\theta e^{-\theta k t}{k\choose k-1}\PP^*(t)^{k} = k\theta U(t)^k$ which cancels with the $- k\theta U(t)^k$ appearing in the first line of the bound. In the last line, we have used the combinatorial identity ${k\choose i} + {k\choose i-1} = {k+1 \choose i}$ for $1 \le i \le k-1$.

Now, using the above bound in \eqref{udiff}, along with the fact that $U(0)=0$ and the induction hypothesis, we obtain for any $t \ge 0$,
\begin{align*}
\mathbb{E}\left(U(t)^k\right) &\le (m+\beta + \theta)\int_0^te^{-\theta s} \sum_{i=0}^{k-1}{k+ 1 \choose i}\mathbb{E}\left(U(s)^{i}\right)ds\\
&\le (m+\beta + \theta)\int_0^te^{-\theta s} \sum_{i=0}^{k-1}{k+ 1 \choose i} D_i \, ds\\
& \le (m+\beta + \theta)D_{k-1} \sum_{i=0}^{k-1}{k+ 1 \choose i}\int_0^te^{-\theta s} ds \\
&\le \frac{m+\beta + \theta}{\theta} D_{k-1}2^{k+1}\\
&= \frac{m+\beta + \theta}{\theta} \cdot \left(\frac{m+\beta + \theta}{\theta}\right)^{k-1} 2^{k(k+1)/2} \cdot 2^{k+1}\\
& = \left(\frac{m+\beta + \theta}{\theta}\right)^{k} 2^{(k+1)(k+2)/2} = D_k.
\end{align*}
This proves \eqref{ubd} for all $\ell \ge 1$, and hence \eqref{lp}.
\end{proof}

\begin{proof}[Proof of Theorem \ref{exponent}]
\eqref{weakcon0} is proved in Theorem \ref{weakconv}.
\eqref{degtail} follows from \cite[Lemma 5.2]{berger2014asymptotic} or \cite[Section 8.4]{van2016random}. 

Recall that $\mathcal{R} = (1-c)(\PP^*(\tau)+1)$, where $\PP^*(\cdot)$ is defined in \eqref{ptdef} and $\tau \sim \operatorname{Exp}(2+ \beta/m)$ is sampled independently of the CTBP $\B^{(m-1+\beta)}(\cdot)$. Thus, it suffices to prove \eqref{prtail} with $\mathcal{R}$ replaced by $\PP^*(\tau)$. 

Write $\alpha := 2 + \beta/m$ and recall that $\theta = 1+(m+\beta)c/m$. Also, from \eqref{martdec} we have that
\begin{equation}\label{pagetau}
\PP^*(\tau) = e^{\theta \tau}\left(M^*(\tau) + \phi^*(\tau)\right).
\end{equation}
From \eqref{phibd}, $\sup_{t \ge 0} \phi^*(t) \le (m+\beta)c/(\theta m) =:\nu$. Hence, for $r \ge 1$,
\begin{align}\label{power1}
\pr{\PP^*(\tau) \ge r} \le \pr{e^{\theta \tau}\left(M^*(\tau) + \nu\right) \ge r} = \int_0^{\infty}\pr{M^*(t) \ge r e^{-\theta t} - \nu} \alpha e^{-\alpha t} dt.
\end{align}
Abbreviating the bound on $\sup_{t \ge 0} \mathbb{E}\left(|M^*(t)|^k\right)$ in \eqref{lp} by $H_k$ and taking any integer $L > \alpha/\theta$, we have that for any $r,t$ such that $r e^{-\theta t} - \nu>0$,
\begin{equation}\label{power2}
\pr{M^*(t) \ge r e^{-\theta t} - \nu} \le \frac{H_L}{\left(r e^{-\theta t} - \nu\right)^L} \wedge 1.
\end{equation}
Note that for any $t \le \frac{1}{\theta} \log \left(\frac{r}{\tilde{C}}\right)$ with $\tilde{C} := \max\{H_L^{1/L} + \nu, 2\nu\}$,
$$
\frac{H_L}{\left(r e^{-\theta t} - \nu\right)^L} \le 1 \ \ \text{ and } \ \ r e^{-\theta t} - \nu \ge r e^{-\theta t}/2.
$$
Hence, using \eqref{power2} in \eqref{power1}, for $r \ge 1$,
\begin{align*}
\pr{\PP^*(\tau) \ge r} &\le \int_0^{\frac{1}{\theta} \log \left(\frac{r}{\tilde{C}}\right)}\frac{H_L}{\left(r e^{-\theta t} - \nu\right)^L} \alpha e^{-\alpha t} dt + \int_{\frac{1}{\theta} \log \left(\frac{r}{\tilde{C}}\right)}^{\infty} \alpha e^{- \alpha t}dt\\
&\le \frac{2^L H_L \alpha}{r^L}\int_0^{\frac{1}{\theta} \log \left(\frac{r}{\tilde{C}}\right)}e^{(\theta L - \alpha)t}dt + \int_{\frac{1}{\theta} \log \left(\frac{r}{\tilde{C}}\right)}^{\infty} \alpha e^{- \alpha t}dt\\
&= \frac{2^L H_L \alpha}{r^L(\theta L - \alpha)} \left(e^{(\theta L - \alpha)\log \left(\frac{r}{\tilde{C}}\right)/\theta} - 1\right) + e^{-\alpha \log \left(\frac{r}{\tilde{C}}\right)/\theta}\\
&\le \frac{2^L H_L \alpha}{r^L(\theta L - \alpha)}\left(\frac{r}{\tilde{C}}\right)^{(\theta L - \alpha)/\theta} + \left(\frac{r}{\tilde{C}}\right)^{-\alpha/\theta}\\
&=\left(\frac{2^L H_L \alpha}{(\theta L - \alpha)\tilde{C}^L} + 1\right)\left(\frac{r}{\tilde{C}}\right)^{-\alpha/\theta}.
\end{align*}
This proves the upper bound in \eqref{prtail}. 

To prove the lower bound, note that by Lemma \ref{lpbdd}, $M^*(\cdot)$ is an $L^2$-bounded martingale and hence, by the martingale convergence theorem \cite[Theorem 11.10]{klenke2013probability},
$$
M^*(t) \xrightarrow{a.s., \ L^2} M^*(\infty) \ \text{ as } \ t \rightarrow \infty
$$
for some random variable $M^*(\infty)$ with finite second moment and zero mean. Moreover, by \eqref{phibd}, $\phi^*(\cdot)$ is non-decreasing and
$
\phi^*(t) \rightarrow \nu \text{ as } t \rightarrow \infty.
$
Hence,
\begin{equation}\label{asconmart}
M^*(t) + \phi^*(t) \xrightarrow{a.s., \ L^2} Z := M^*(\infty) + \nu \ \text{ as } \ t \rightarrow \infty,
\end{equation}
where $\mathbb{E}(Z) = \nu>0$. Hence, there exist positive $t_0, \eta_1,\eta_2$ such that for all $t \ge t_0$,
$$
\pr{M^*(t) + \phi^*(t) > \eta_1} \ge \frac{1}{2}\pr{Z> \eta_1} \ge \eta_2.
$$
Using this and \eqref{pagetau}, we obtain for $r \ge \eta_1 e^{\theta t_0}$,
\begin{align*}
\pr{\PP^*(\tau) \ge r} &\ge \int_{\frac{1}{\theta} \log \left(\frac{r}{\eta_1}\right)}^{\infty}\pr{M^*(t) + \phi^*(t)  \ge r e^{-\theta t}} \alpha e^{-\alpha t} dt\\
&\ge \int_{\frac{1}{\theta} \log \left(\frac{r}{\eta_1}\right)}^{\infty}\pr{M^*(t) + \phi^*(t) > \eta_1} \alpha e^{-\alpha t} dt\\
&\ge \int_{\frac{1}{\theta} \log \left(\frac{r}{\eta_1}\right)}^{\infty}\eta_2 \alpha e^{-\alpha t} dt = \eta_2 e^{-\frac{\alpha}{\theta} \log \left(\frac{r}{\eta_1}\right)} = \eta_2 \left(\frac{r}{\eta_1}\right)^{-\alpha/\theta}.
\end{align*}
This proves the lower bound in \eqref{prtail} and completes the proof of the theorem. 
\end{proof}

\section{Proofs: PageRank asymptotics for the root when $m=1$}\label{rootsec}
In this section, we will prove Theorem \ref{rootpr}.

Recall the directed random graph process $\{\G_n: n \ge 1\}$ defined in Section \ref{modeldef} with $\beta \ge 0$ and $m=1$. Also recall the CTBP $\{\B^{(\beta)}(t) : t \ge 0\}$ defined after Definition \ref{ctbp}. Define the following stopping times: $T_0 =0$ and for $n \ge 1$,
$
T_n := \inf\{t \ge 0 : |\B^{(\beta)}(t)| = n\}.
$

An important connection between CTBP and the discrete random tree sequence $\{\G_n\}_{n \ge 1}$ with $m =1$ is given by the following result which is a consequence of the properties of the exponential distribution (and is the starting point of the Athreya-Karlin embedding \cite{athreya1968}). 
\begin{lemma}\label{lem:ctb-embedding}
 Viewed as a sequence of growing random labelled directed rooted trees,
 $$
\{\B^{(\beta)}(T_n) : n \ge 0\} \equald \{\G_n : n \ge 1\}.
 $$
\end{lemma}
Let $\rho(\cdot)$ denote the intensity measure of a $\beta$-Yule process. Define its (formal) Laplace transform $\hat{\rho}: (0,\infty) \rightarrow (0, \infty]$ by
$$
\hat{\rho}(\lambda) := \int_0^{\infty} e^{-\lambda t} \rho(t) dt = \sum_{n=1}^{\infty} \prod_{i=0}^{n-1}\frac{i+\beta}{\lambda + i + \beta}.
$$
Explicit computation shows that $\hat{\rho}(\lambda) = (1+\beta)/(\lambda-1), \ \lambda >1$ (see \cite[Section 4.2]{rudas-2}, but note that the formula for $\hat{\rho}$ there is slightly different since the attachment probabilities are proportional to $\beta$ plus the \emph{in-degree} of vertices, while here they are proportional to $\beta$ plus the \emph{total degree}).  Hence, the equation $\hat{\rho}(\lambda)=1$ has a unique positive root given by $\lambda^* = 2+ \beta$.
$\lambda^*$ is called the \emph{Malthusian rate of growth} as it captures the asymptotic exponential growth rate of the branching process population size $|\B^{(\beta)}(\cdot)|$ (see \eqref{growas}). We refer the reader to \cite{rudas-2} for a detailed treatment on the connection between random tree processes and CTBPs.

The analysis of the CTBP using Malthusian rates goes back to \cite{jagers-nerman-1,jagers-nerman-2,nerman1981convergence}. We will use the following fact \cite[Theorem 5.4]{nerman1981convergence}:
\begin{equation}\label{growas}
e^{-(2+\beta) t} |\B^{(\beta)}(t)| \xrightarrow{a.s.} \Phi \ \text{ as } \ t \rightarrow \infty
\end{equation}
for some positive, almost surely finite random variable $\Phi$.

\begin{proof}[Proof of Theorem \ref{rootpr}]
By Lemma \ref{lem:ctb-embedding},
$$
\{ R_1(n) : n \ge 1\} \equald \{(1-c)(1 + \PP^*(T_n)) : n \ge 1\},
$$
where $\PP^*(\cdot)$ is defined in \eqref{ptdef}. Thus, it suffices to show \eqref{rootas} with $\PP^*(T_n)$ in place of $R_1(n)$.

By \eqref{martdec}, recalling $\theta = 1+(1+\beta)c$ when $m=1$,
\begin{equation}\label{r1}
\PP^*(t) = e^{\theta t}\left(M^*(t) + \phi^*(t)\right), \ t \ge 0.
\end{equation}
Recall from \eqref{asconmart} that
\begin{equation}\label{r2}
M^*(t) + \phi^*(t) \xrightarrow{a.s.} Z := M^*(\infty) + \nu \ \text{ as } \ t \rightarrow \infty,
\end{equation}
where $\nu = (1+\beta)c/\theta$ and $Z$ is a random variable with finite second moment and $\mathbb{E}(Z) = \nu>0$. By \eqref{growas},
\begin{equation}\label{r3}
e^{-(2+\beta)T_n} n = e^{-(2+\beta)T_n}  |\B^{(\beta)}(T_n)| \xrightarrow{a.s.} \Phi \ \text{ as } \ n \rightarrow \infty.
\end{equation}
Hence, by \eqref{r2} and \eqref{r3},
\begin{align*}
n^{-\theta/(2+\beta)}\PP^*(T_n) &= n^{-\theta/(2+\beta)}e^{\theta T_n}\left(e^{-\theta T_n}\PP^*(T_n)\right)\\
&= \left(e^{-(2+\beta)T_n} n\right)^{-\theta/(2+\beta)}\left(e^{-\theta T_n}\PP^*(T_n)\right)\\
&\xrightarrow{a.s} \Phi^{-\theta/(2+\beta)} Z, \ \text{ as } \ n \rightarrow \infty.
\end{align*}
This proves the theorem.
\end{proof}

\bibliographystyle{imsart-number}
\bibliography{scaling}
\end{document}